\documentclass{article}
\usepackage[utf8]{inputenc}
\usepackage{amsmath,amsthm}
\usepackage{amssymb}
\usepackage{xspace}
\usepackage{xcolor}
\usepackage{url}
\usepackage{graphicx}
\usepackage{complexity}
\usepackage{caption}
\usepackage{subcaption}
\usepackage{multirow} 
\usepackage{lscape} 
\usepackage[affil-it]{authblk} 
\usepackage{hyperref}
\usepackage{enumitem} 

\newcommand{\N}{{\mathbb{N}\setminus\{0\}}}

\newcommand{\blue}[1]{#1} 

\newcommand{\twalk}[1]{temporal $#1$-walk}
\newcommand{\tpath}[1]{temporal $#1$-path}
\newcommand{\vcut}[1]{temporal $#1$-cut}

\newtheorem{theorem}{Theorem}
\newtheorem{problem}[theorem]{Problem}
\newtheorem{lemma}[theorem]{Lemma}
\newtheorem{proposition}[theorem]{Proposition}
\newtheorem{claim}{Claim}[theorem]

\title{Mengerian graphs: characterization and recognition\thanks{Funded by CNPq GD 141498/2021-8, CNPq Produtividade 303803/2020-7 and FUNCAP PS1-00186-00155.01.00/21.\\E-mail: \texttt{allen.ibiapina@alu.ufc.br,anasilva@mat.ufc.br}.}}

\author{Allen Ibiapina \and Ana Silva\thanks{This work was developed during this author's visit to Dipartimento di Sistemi, Informatica, Applicazioni - Universit\`{a} degli Studi di Firenze, Italy.}}
\affil{ParGO Group - Parallelism, Graphs and Optimization\\ Departamento de Matemática - Centro de Ciências\\ Universidade Federal do Ceará, Brazil}

\date{}

\begin{document}

	\maketitle
	
	\begin{abstract}
	A temporal graph ${\cal G}$ is a graph that changes with time. More specifically, it is a pair $(G, \lambda)$ where $G$ is a graph and $\lambda$ is a function on the edges of $G$ that describes when each edge $e\in E(G)$ is active. Given vertices $s,t\in V(G)$, a \tpath{s,t} is a path in $G$ that traverses edges in non-decreasing time; and if $s,t$ are non-adjacent, then a \vcut{s,t} is a subset $S\subseteq V(G)\setminus\{s,t\}$ whose removal destroys all \tpath{s,t}s.
    
    It is known that Menger's Theorem does not hold on this context, i.e., that the maximum number of internally vertex disjoint \tpath{s,t}s is not necessarily equal to the minimum size of a \vcut{s,t}. 
    In a seminal paper, Kempe, Kleinberg and Kumar (STOC'2000) defined a graph $G$ to be Mengerian if equality holds on $(G,\lambda)$ for every function $\lambda$. They then proved that, if each edge is allowed to be active only once in $(G,\lambda)$, then $G$ is Mengerian if and only if $G$ has no gem as topological minor. In this paper, we generalize their result by allowing edges to be active more than once, giving a characterization also in terms of forbidden structures. We additionally provide a polynomial time recognition algorithm.
	\end{abstract}
	
	\section{Introduction}
	Temporal graphs have been the subject of a lot of interest in recent years (see e.g. the surveys ~\cite{Holme.15,LVM.18,M.16,Netal.13}). They appear under many distinct names (temporal networks~\cite{Holme.15}, edge-scheduled networks~\cite{B.96}, dynamic networks~\cite{XFJ.03}, time-varying graphs~\cite{CFQS.12}, stream graphs, link streams~\cite{LVM.18}, etc), but with very little  (if any) distinction between the various models. 
Here, we favor the name temporal graphs.

An example where one can apply a temporal graph is the modeling of proximity of people within a region, with each vertex representing a person, and two people being linked by an edge at a given moment if they are close to each other. This and similar ideas have been used also to model animal proximity networks, human communication, collaboration networks, travel and transportation networks, etc. We refer the reader to~\cite{Holme.15} for a plethora of applications, but just to cite a simple (and trendy) example where these structures could be used, imagine one wants to track, in the proximity network previously described, the spreading of a contagious disease. In this case it is important to be able to detect whether there was an indirect contact between two people, and of course the contact is only relevant if it occurred after one of these people got sick (see e.g.~\cite{EK.18}). 
Therefore, when studying temporal graphs, it makes sense to be concerned with paths between vertices that respect the flow of time; these are called \emph{temporal paths} and are essentially different from the traditional concept in static graphs, allowing even for multiple definitions of minimality (see~e.g.~\cite{CHMZ.19,MMS.19,XFJ.03}).

Temporal graphs are being used in practice in a variety of fields for many decades now, with the first appearances dating back to the 1980's, but only recently there seems to be a bigger effort to understand these structures from a more theoretical point of view. An issue that is often raised is whether results on the static version of a certain problem are inherited. 
This is not always the case, as has been shown for the temporal versions of connected components~\cite{BF.03}, Menger's Theorem~\cite{B.96,KKK.00,ILMS.arxiv}, Edmonds' Theorem~\cite{CLMS.20,KKK.00} and Euler's Theorem~\cite{BM.21,MS.21}. 

In particular, in the seminal paper~\cite{KKK.00}, Kempe, Kleinberg and Kumar define and characterize a Mengerian graph as being a graph\footnote{We adopt the definition of~\cite{West.book}, where a graph can have multiple edges incident to the same pair of vertices. This is sometimes called multigraph.} for which all assignment of activity times for its edges produces a temporal graph on which the temporal version of Menger's Theorem holds. However, their definition does not allow for a certain edge to be active multiple times, which is the case in most of the recent studies. In this paper, we fill this gap, giving a complete characterization of Mengerian graphs. Like the characterization given in~\cite{KKK.00}, ours is in terms of forbidden structures. We also provide a polynomial-time  recognition algorithm for these graphs.

This work has been presented in~\cite{IS.21_FCT}, and here we provide complete proofs of our results.

\paragraph{Our Results.}

Given a graph $G$, and a time-function $\lambda:E(G)\rightarrow \N$, we call the pair $(G,\lambda)$ a \emph{temporal graph}. Also, given vertices $s,t\in V(G)$, a \emph{\tpath{s,t}} in $(G,\lambda)$ is a path between $s$ and $t$ in $G$ such that if edge $e$ appears after $e'$ in the path, then $\lambda(e') \le \lambda(e)$. We say that two such paths are \emph{vertex disjoint} if their internal vertices are distinct. Also, a subset $S\subseteq V(G)\setminus\{s,t\}$ is a \emph{\vcut{s,t}} if there is no \tpath{s,t} in $(G-S,\lambda')$, where $\lambda'$ is equal to $\lambda$ restricted to $E(G-S)$. From now on, we make an abuse of language and write simply $(G-S,\lambda)$. 
The analogous notions on static graphs, and on edges, can be naturally defined and the following is a classic result in Graph Theory.

\begin{theorem}[Menger~\cite{Menger.27}]
Let $G$ be a graph, and $s,t\in V(G)$ be such that $st\notin E(G)$. Then, the maximum number of vertex (edge) disjoint $s,t$-paths in $G$ is equal to the minimum size of a vertex (edge) cut in $G$.
\end{theorem}

\begin{figure}[h]
		\centering
		\begin{subfigure}[b]{0.4\textwidth}
			\includegraphics[width=5cm]{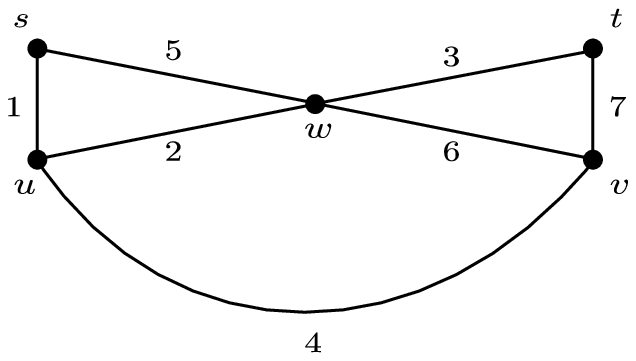}
			\caption{Simple graph.}\label{fig:KKKexample}
		\end{subfigure}
		\hfill
		\begin{subfigure}[b]{0.4\textwidth}
			\includegraphics[width=5cm]{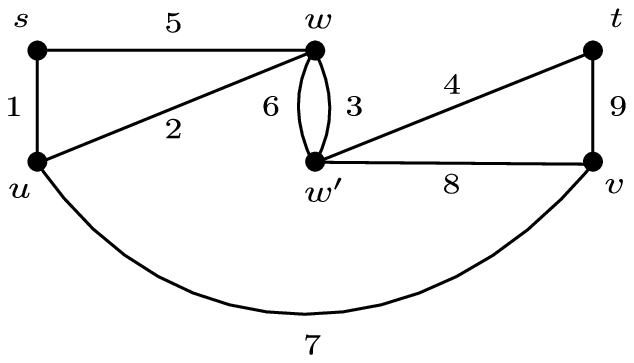}
			\caption{General graph.}\label{fig:Generalexample}
		\end{subfigure}
		\caption{Examples of graphs where the temporal version of Menger's Theorem applied to vertices does not hold.}
	\end{figure}

In~\cite{B.96}, the author already pointed out that this might not be the case for temporal paths. Here, we present the example given later in~\cite{KKK.00}; see Figure~\ref{fig:KKKexample}. Observe that the only \tpath{s,t} using $sw$  also uses $v$, and since $\{w,v\}$ separates $s$ from $t$, we get that there are no two vertex disjoint \tpath{s,t}s. At the same time, no single vertex in $\{u,v,w\}$ breaks all \tpath{s,t}s, i.e., there is no \vcut{s,t} of size~1.

Inspired by this, in~\cite{KKK.00} the authors define the \emph{Mengerian graphs}, which are those graphs $G$ for which the temporal version of Menger's Theorem always holds, i.e., for every time-function $\lambda$, we get that the maximum number of \tpath{s,t}s in $(G,\lambda)$ is equal to the minimum size of a \vcut{s,t} in $(G,\lambda)$, for every pair $s,t\in V(G)$ such that $st\notin E(G)$. They then give a characterization for when the graph is simple: a simple graph $G$ is Mengerian if and only if $G$ does not have the graph depicted in Figure~\ref{fig:KKKexample} (also known as the gem) as a topological minor. Despite being a very nice result, it does not allow for an edge of a simple graph to be active more than once, which is generally the case in practice. To see this, observe the graph depicted in Figure~\ref{fig:Generalexample}. This graph does not have the gem as a topological minor, and is non-Mengerian since the given time-function is such that there is only one \tpath{s,t}, while the minimum size of a \vcut{s,t} is two (this can be seen by similar arguments as the ones applied to Figure~\ref{fig:KKKexample}).

In this paper, we allow $G$ to be a general graph, i.e., it can contain multiple edges incident to the same pair of vertices (this is also sometimes called multigraph). Note that the assignment of time labels to such a graph can be seen as allowing for the edges of a simple graph to be active more than once in time. We prove the characterization below. The formal definition of an m-topological minor is given in Section~\ref{sec:defs}, but for now it suffices to say that, when subdividing an edge $e$, the multiplicity of the obtained edges is the same as $e$.

\begin{figure}
		\centering
		\begin{subfigure}[b]{0.4\textwidth}
			\centering
			\includegraphics[width=\textwidth]{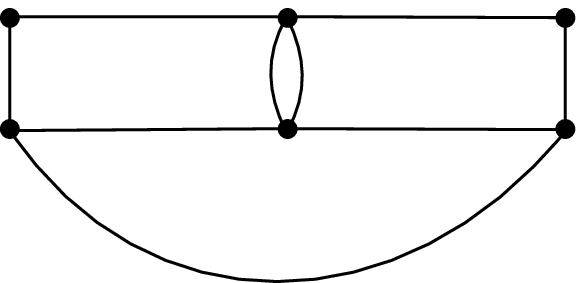}
			\caption{${\cal F}_1$}
		\end{subfigure}
		\hfill
		\begin{subfigure}[b]{0.4\textwidth}
			\centering
			\includegraphics[width=\textwidth]{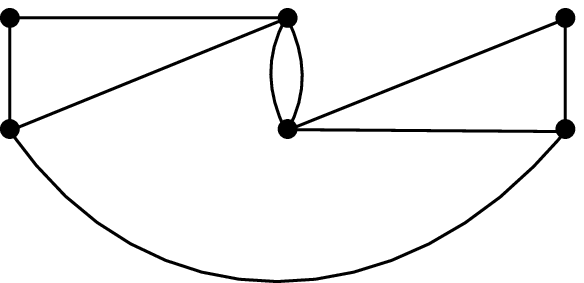}
			\caption{${\cal F}_2$}
		\end{subfigure}
		\hfill
		\begin{subfigure}[b]{0.4\textwidth}
			\centering
			\includegraphics[width=\textwidth]{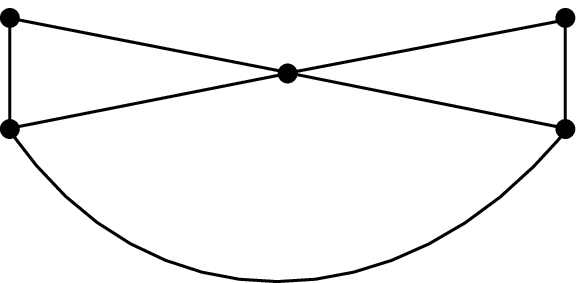}
			\caption{${\cal F}_3$}
		\end{subfigure}
		\caption{Forbidden m-topological minors.}
		\label{fig:ForbMinors}
	\end{figure}

\begin{theorem}\label{main}
Let $G$ be any graph. Then, $G$ is a Mengerian graph if and only if $G$ does not have one of the graphs in Figure~\ref{fig:ForbMinors} as an m-topological minor.
\end{theorem}

We also provide a polynomial time recognition algorithm to decide whether a given graph $G$ is Mengerian. 
\begin{theorem}\label{thm:recog}
Let $G$ be a graph on $n$ vertices and $m$ edges. Then, one can decide whether $G$ is a Mengerian graph in time $O(mn^3)$.
\end{theorem}

\paragraph{Related works and remarks.} 
In~\cite{B.96}, Berman proves that the temporal version of Menger's Theorem applied to temporal edges always holds. It should be noted that Berman uses the same model as us, i.e., an edge has always a single appearance and multiple edges are allowed. If instead the paths are not allowed to share distinct edges that have the same endpoints (i.e., they are multiedge disjoint), then Menger's Theorem also does not necessarily hold. To see this, observe again Figure~\ref{fig:Generalexample}. Note that: every \tpath{s,t} passing by $sw$ must use $ww'$ (at time 6) and $vt$; and every \tpath{s,t} passing by $su$ uses either $ww'$ (at any time) or $vt$. Therefore there are no two multiedge disjoint \tpath{s,t}s. Because there is no $xy\in E(U(G))$ such that the removal of every edge with endpoints $xy$ breaks all the \tpath{s,t}s, we get that the multiedge version of Menger's Theorem also does not hold.  Then it is natural to ask whether there is a characterization of Mengerian graphs on the multiedge context.
	
	\begin{problem}
	Let ${\cal G}$ contain every $G$ for which the multiedge version of Menger's Theorem holds on $(G,\lambda)$ for every $\lambda$. Can ${\cal G}$ be characterized in terms of forbidden structures?
	\end{problem}

In~\cite{MMS.19}, Mertzios, Michail and Spirakis give an alternative formulation of Menger's Theorem that holds on the temporal context. 
There, they define the notion of \emph{out-disjointness}, where two paths are disjoint if they do not share the same departure time for a given vertex, and the notion of \emph{node departure time cut}, where one removes a time label from the possible departure times of a vertex. They then prove that the maximum number of out-disjoint paths between $s,t$ is equal to the minimum size of a node departure time cut. 
In~\cite{B.96}, Berman proves that deciding whether there are $k$ vertex disjoint \tpath{s,t}s is $\NP$-complete. This was improved for fixed $k=2$ in~\cite{KKK.00}, where Kempe, Kleinberg and Kumar also prove that deciding whether there is a \vcut{s,t} of size at most $k$ is $\NP$-complete, for given $k$. Observe that the latter problem can be easily solved in time $O(|V(G)|^k)$, which raises the question about whether it is $\FPT$ when parameterized by $k$. This is answered negatively in~\cite{ZFMN.20}, where Zschoche, Fluschnik, Molter and Niedermeier prove that this is $\W[1]$-hard. 
In~\cite{FMNRZ.20}, Fluschnik, Molter, Niedermeier, Renken and Zschoche further investigate the cut problem, giving a number of hardness results (e.g. that the problem is hard even if $G$ is a line graph), as well as some positive ones (e.g., that the problem is polynomial-time solvable when $G$ has bounded treewidth). 
Finally, in~\cite{ILMS.arxiv}, Ibiapina, Lopes, Marino and Silva investigate the temporal vertex disjoint version of Menger's Theorem and, among other results, prove that the related paths and walks problems differ and that Menger's Theorem applied to temporal vertex disjoint paths holds if and only if the maximum number of such paths is equal to~1.

Our paper is organized as follows. 
We give the formal definitions in Section~\ref{sec:defs}, and present the proofs in Section~\ref{sec:necessity} (necessity of Theorem~\ref{main}), Section~\ref{sec:sufficiency} (sufficiency of Theorem~\ref{main}), and Section~\ref{sec:recog} (proof of Theorem~\ref{thm:recog}).

\section{Basic Definitions and Terminology}\label{sec:defs}
	
	Given a natural number $n$, we denote by $[n]$ the set $\{1,2,\ldots,n\}$. A \emph{graph} is an ordered triple $G= (V,E,f)$ where $V$ is a finite non-empty set (called \emph{set of vertices}), $E$ is a finite set (called the \emph{set of edges}), and $f$ is a function that associates to each $e \in E$ an unordered pair $xy$ of vertices of $V$. The pair $xy$ is called the \emph{endpoints of $e$}, and if a pair $xy$ has at least one related edge, then $xy$ is also called \emph{multiedge}. For simplicity, in what follows we omit the function $f$ and simply refer to the endpoints of each edge instead. The \emph{multiplicity} of $xy$ is the number of edges which are associated with this pair. If the graph is denoted by $G$, we also use $V(G),E(G)$ to denote its vertex and edge sets, respectively. A graph $H$ such that $V(H)\subseteq V(G)$ and $E(H)\subseteq V(G)$ is called a \emph{subgraph of $G$}, and we write $H\subseteq G$. The subgraph obtained by removing edges in such way that all multiedges has multiplicity one is called \emph{underlying simple graph} of $G$ and is denoted by $U(G)$.

	The \emph{degree} of a vertex $v \in V(G)$ is the number of edges having $v$ as endpoint, and it is denoted by $d_{G}(v)$. The \emph{maximum degree} of $G$ is then the maximum among the degree of vertices of $G$; it is denoted by $\Delta(G)$. Given two vertices $s,t \in V(G)$, an \emph{$s,t$-walk} $P$ in $G$ is a sequence $(s=u_{1},e_{1},\ldots,e_{q-1},u_{q} = t)$ that alternates vertices and edges and is such that $e_{i}$ is an edge with endpoints $u_{i}u_{i+1}$, for every $i \in [q-1]$. If additionally no vertex appears more than once in $P$ then $P$ is called an \emph{$s,t$-path}; and if $u_1=u_q$ and this is the only vertex appearing more than once in $P$, then $P$ is called a \emph{cycle}. If $G$ is a simple graph, then we omit the edges in the sequence $P$. 
	We denote by $V(P)$ and $E(P)$ the sets of vertices and edges appearing in $P$, respectively. 
	
\begin{figure}[h]
		\centering
		\includegraphics[width = 8cm]{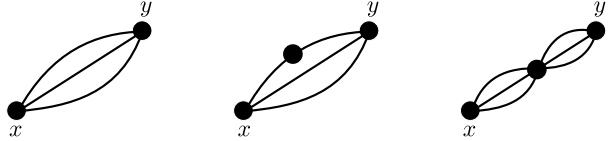}
		\caption{From left to right: the multiedge $xy$, the subdivision of an edge with endpoints $xy$, and the m-subdivision of $xy$.}
		\label{fig:defmmenor}
	\end{figure}

	Given a simple graph $G$ and vertices $Z\subseteq  V(G)$, the \emph{identification of $Z$} is the graph obtained from $G - Z$ by adding a new vertex $z$ and, for every edge $e$ with endpoints $z'u$ where $z' \in Z$ and $u\notin Z$, add an edge $e'$ with endpoints $zu$. The graph $G'$ obtained from $G$ by a \emph{subdivision} of an edge $e$ with endpoints $uv$ is the graph having $V(G)\cup \{z_e\}$ as vertex set, and $E(G-e)\cup\{e', e''\}$ as edge set, where $e'$ has endpoints $uz_e$ and $e''$ has endpoints $z_ev$. Finally, the graph obtained from $G$ by an \emph{m-subidivision} of a multiedge $xy$ is the graph obtained by subdividing all the edges with endpoints $xy$ and then identifying the new vertices. Observe Figure~\ref{fig:defmmenor} for an illustration of these definitions. Given a graph $H$, if $G$ has a subgraph that can be obtained from m-subdivisions of $H$, then we say that $H$ is an \emph{m-topological minor} of $G$ and write $H \preceq_m G$. And if $G$ has a subgraph that can be obtained from subdivisions of $H$, then we say that $H$ is a \emph{topological minor} of $G$ and write $H \preceq G$

	A \emph{chain} in $G$ is a path $(u_1,\cdots,u_q)$ in $U(G)$, $q\ge 2$, such that the multiplicity of $u_{i}u_{i+1}$ in $G$ is at least~two, for each $i\in [q-1]$. Observe that each m-subdivision $H$ of ${\cal F}_1$ or ${\cal F}_2$ contains exactly one chain $P$; we call $P$ the \emph{chain of $H$}.
	
	A \emph{temporal graph} is a pair $(G,\lambda)$ where $G$ is a graph and $\lambda$ is a function that assigns to each edge $e$ a value in $\N$. We say that $e$ is \emph{active in time $\lambda(e)$}, and call $\lambda$ the \emph{time-function} of $(G,\lambda)$.  
	Given $s,t \in V(G)$, a \emph{\twalk{s,t}} $P$ is a walk $(s,e_{1},v_{2},e_{2},\ldots,e_{k-1},t)$ in $G$ such that $\lambda(e_{1}) \leq \ldots \leq \lambda(e_{k-1})$; if $P$ is also a path in $G$, then it is called a \emph{\tpath{s,t}}. 
	To emphasize the time appearances of such edges, we sometimes denote $P$ by $(s,\lambda(e_{1}),v_{2},\lambda(e_{2}),\ldots,\lambda(e_{k-1}),t)$. For a pair of vertices $uw$ we denote by $\lambda(uw)$ the set $\{\lambda(e)~\colon~ \text{e has endpoints }uw\}$; if such set is unitary we abuse the notation by saying that $\lambda(uw)$ is the value in such set. 
	
	For a subset $S$ of vertices or edges, we denote by $(G-S,\lambda)$ the temporal graph $(G-S,\lambda')$, where $\lambda'$ is the restriction of $\lambda$ to $E(G-S)$. Similarly, for a subgraph $H\subseteq G$, the temporal graph $(H,\lambda)$ denotes $(H,\lambda')$, where $\lambda'$ is the restriction of $\lambda$ to $E(H)$.
	
	Let $Z\subseteq V(G)$ be any subset of vertices. The \emph{identification of $Z$ in $(G,\lambda)$} is the temporal graph $(G^*,\lambda^*)$, where $G^*$ is the graph obtained from $G$ by the identification of $Z$ into $z$, and $\lambda^*$ is such that $\lambda^*(e) = \lambda(e)$ if $e \in E(G)$ has no endpoints in $Z$, and $\lambda^*(uz)=\bigcup_{x \in Z} \lambda(ux)$ for every $u \in N(Z)\setminus Z$.
	
	
	Given two \twalk{s,t}s $P$ and $Q$, we say that $P$ and $Q$ are \emph{disjoint} if they do not share any internal vertex, i.e., if $V(P)\cap V(Q)=\{s,t\}$. Also, if $s,t$ are non-adjacent, then a \emph{\vcut{s,t}} is a set $S \subseteq V(G)\setminus \{s,t\}$ such that there are no temporal $s,t$-walks in $(G-S,\lambda)$. 
	Throughout the text, whenever we talk about a \vcut{s,t} we implicitly assume that it is well defined (in other words, that $s,t$ are non-adjacent).
	Denote by $p_{G,\lambda}(s,t)$ the maximum number of vertex disjoint \twalk{s,t}s and by $c_{G,\lambda}(s,t)$ the size of a minimum \vcut{s,t}. If $(G,\lambda)$ is clear from the context, we omit it from the subscript. The following proposition can be easily seen to hold and tells us that these definitions can be also made in terms of paths instead of walks. 
	 
	\begin{proposition}\label{prop:walk}
	Let $(G,\lambda)$ be a temporal graph and $s,t \in V(G)$. Then every \twalk{s,t} contains a \tpath{s,t}. Additionally, if $p^*(s,t)$ and $c^*(s,t)$ are defined as above, but concerning temporal paths instead of temporal walks, then $p^*(s,t) = p(s,t)$ and $c^*(s,t)=c(s,t)$. 
    \end{proposition}
	
	

	In the remainder of the text, we sometimes refer to a \twalk{s,t} as being a \tpath{s,t}; it means that we are implicitly dealing with the path contained within the related walk.  
	
	 


	 Throughout the text, we make concatenations of paths and of temporal paths and for this we need some further notation. If $P$ is a (temporal) path and $u,v\in V(P)$, then the (temporal) $u,v$-path contained in $P$ is denoted by $uPv$. If $P$ is a $u,v$-path and $Q$ is a $v,w$-path, then we denote the $u,w$-walk obtained from the concatenation of $P$ and $Q$ by $PQ$. If additionally $P$ and $Q$ are temporal paths, and $P$ finishes at most at the time when $Q$ starts, then $PQ$ denotes the temporal $u,w$-walk obtained from the concatenation of $P$ and $Q$.

	
	\section{Proof of necessity of Theorem~\ref{main}}\label{sec:necessity}

	\begin{figure}
		\centering
		\begin{subfigure}[b]{0.4\textwidth}
			\centering
			\includegraphics[width=\textwidth]{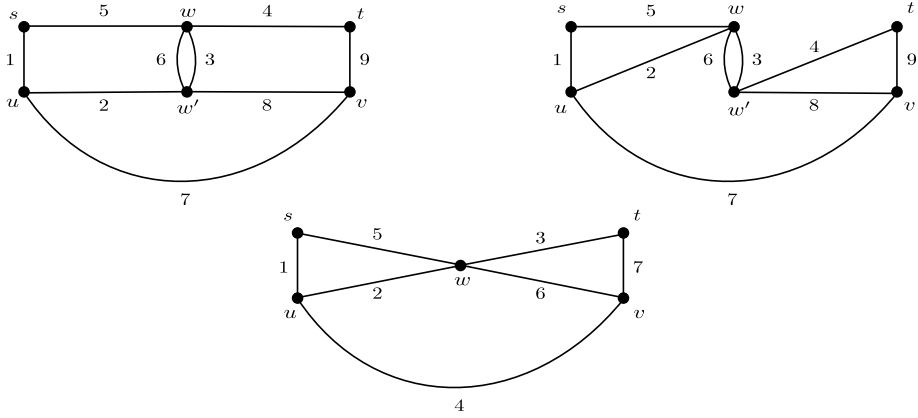}
		\end{subfigure}
		\hfill
		\begin{subfigure}[b]{0.4\textwidth}
			\centering
			\includegraphics[width=\textwidth]{F_3-rotulado.eps}
		\end{subfigure}
		\hfill
		\begin{subfigure}[b]{0.4\textwidth}
			\centering
			\includegraphics[width=\textwidth]{F_4-rotulado.eps}
		\end{subfigure}
		\caption{Forbidden m-topological minors together with time-functions that show that these are non-Mengerian graphs.}
		\label{fig:badlabeling}
	\end{figure}
	
	In this section, we prove that if $G$ has an m-subdivision of $\mathcal{F}_1$, $\mathcal{F}_2$ or $\mathcal{F}_3$, then $G$ is non-Mengerian. We start by proving the following lemma.
	
	\begin{lemma}\label{lem:nonmeng}
	    Graphs $\mathcal{F}_1$, $\mathcal{F}_2$ and $\mathcal{F}_3$ are non-Mengerian.
	\end{lemma}
	
	\begin{proof}
	    We use the time-functions represented in Figure~\ref{fig:badlabeling}. It was proved in~\cite{KKK.00} that $\mathcal{F}_3$ is non-Mengerian. We show that $\mathcal{F}_1$ and $\mathcal{F}_2$ are non-Mengerian similarly. To show that $c(s,t)\geq 2$ in $\mathcal{F}_1$ and $\mathcal{F}_2$, we show that for each vertex $x$ different from $s,t$ there is a \tpath{s,t} not passing by $x$. For $w$ and $w'$ we take the temporal path $(s,1,u,7,v,9,t)$. A \tpath{s,t} not passing by $u$ is $(s,5,w,6,w',8,v,9,t)$. Finally, a \tpath{s,t} not passing by $v$ in $\mathcal{F}_1$ is $(s,1,u,2,w',3,w,4,t)$, while such a path in $\mathcal{F}_2$ is $(s,1,u,2,w,3,w',4,t)$. 
	    Now we show that $p(s,t) = 1$ in both $\mathcal{F}_1$ and $\mathcal{F}_2$. Suppose otherwise and observe that $\{v,w\}$ is a \vcut{s,t}. Note also that any 2 disjoint \tpath{s,t}s must include a path that starts with the edge labelled with~5; let $P$ be such path. Notice, however, that the subgraph formed by the edges active at time at least~5 is a tree, and hence contains exactly one $s,t$-path. It thus follows that $P$ must be equal to the path contained in such tree, in which case $\{v,w\}\subseteq V(P)$. As $\{v,w\}$ is a \vcut{s,t}, any other \tpath{s,t} must intersect $P$, a contradiction. 
	\end{proof}

	The following lemma shows us that any m-subdivision of ${\cal F}_i$ is also non-Mengerian, for every $i\in \{1,2,3\}$.
	
	\begin{lemma}\label{lem:msubdiv}
		If $G$ is non-Mengerian, then an m-subdivision of $G$ is also non-Mengerian. 
	\end{lemma}
	\begin{proof}
		Let $G$ be a non-Mengerian graph and consider $\lambda \in E(G) \to \N$ and $s,t \in V(G)$ to be such that $p(s,t) < c(s,t)$. Also, suppose that $H$ is obtained from $G$ by m-subdividing a multiedge, say $xy$. We construct a function $\lambda'$ from $\lambda$ that proves that $H$ is also non-Mengerian. 
		
    		Let $D\subseteq E(G)$ be the set of edges of $G$ with endpoints $xy$, and denote by $v_{xy}$ the vertex of $H$ created by the m-subdivision of $xy$.  Moreover, denote by $D_x$ and $D_y$ the sets of edges of $H$ with endpoints $xv_{xy}$ and $v_{xy}y$, respectively. Finally, define $\lambda'$ to be such that  $\lambda'(e) = \lambda(e)$, for every $e \in E(G)\setminus D$, and $\lambda'(D_x) = \lambda'(D_y) = \lambda(D)$. We show that $c_{G,\lambda}(s,t) = c_{H,\lambda'}(s,t)$ and $p_{G,\lambda}(s,t) = p_{H,\lambda'}(s,t)$, which finishes our proof. 
    		
    		Given a set of disjoint \tpath{s,t}s in $(G,\lambda)$, if some of these paths, say $P$, uses the edge $xy$, then in $(H,\lambda')$ we can substitute such edge by an edge in $D_x$ and another in $D_y$ with the same time to obtain a temporal path $P'$ such that $V(P')=V(P)\cup \{v_{xy}\}$. This gives us a set of disjoint \tpath{s,t}s in $(H,\lambda')$. In the other direction, if it is given a set of disjoint \tpath{s,t}s in $(H,\lambda')$, if some of them uses the vertex $v_{xy}$, let $f_j$ be the edge used in $D_j$ for $j \in \{x,y\}$. Suppose without loss of generality that $\lambda(f_x)\geq \lambda(f_y)$. Then we substitute both edges incident to $v_{xy}$ by an edge in $D$ appearing at time $\lambda(f_x)$. This implies that $p_{G,\lambda}(s,t) = p_{H,\lambda'}(s,t)$. \blue{To see that $c_{G,\lambda}(s,t) = c_{H,\lambda'}(s,t)$ one can just additionally notice that if a \vcut{s,t} in $(H,\lambda')$ contains $v_{xy}$, then $(S\setminus v_{xy})\cup \{x\}$ is a \vcut{s,t} in $(G,\lambda)$. Also, that any \vcut{s,t} in $(G,\lambda)$ is also a \vcut{s,t} in $(H,\lambda')$}.
	\end{proof}

	Recall that we want to prove necessity of Theorem~\ref{main}, i.e., that if $G$ has $\mathcal{F}_1,\mathcal{F}_2$ or $\mathcal{F}_3$ as an m-topological minor, then $G$ is non-Mengerian. Observe that if $G$ has $\mathcal{F}_1,\mathcal{F}_2$ or $\mathcal{F}_3$ as an m-topological minor, then either $G$ is itself an m-subdivision of one of these graphs, or it contains a subgraph $H$ that is an m-subdivision of one of these graphs. 
	Since Lemmas~\ref{lem:nonmeng} and~\ref{lem:msubdiv} already give us that any m-subdivision of $\mathcal{F}_1,\mathcal{F}_2$ or $\mathcal{F}_3$ is non-Mengerian, we finish our proof with the next lemma.

	\begin{lemma}\label{lemma:subgraphClosed}
		$G$ is Mengerian if and only if $H$ is Mengerian, for every $H\subseteq G$.
	\end{lemma}
	\begin{proof}
		To prove necessity, suppose that $H\subseteq G$ is non-Mengerian, and let $s,t,\lambda$ be such that $p_{H,\lambda}(s,t) < c_{H,\lambda}(s,t)$. Consider the time-function  $\lambda'$ for $E(G)$ defined as follows. 
		\[\lambda(e) = \left\{\begin{array}{ll}
			\lambda(e)+1 & \mbox{, for every } e\in E(H),\\
			1 & \mbox{, for every $e\in E(G)\setminus E(H)$ with endpoints $yt$, and}\\
			\max\lambda(E(H))+2 & \mbox{, otherwise.}
		\end{array}\right.\]

		Because $H\subseteq G$ and $\lambda\subseteq \lambda'$, note that we get $p_{H,\lambda}(s,t)\le p_{G,\lambda'}(s,t)$ and $c_{H,\lambda}(s,t)\le c_{G,\lambda'}(s,t)$. Therefore it suffices to prove $p_{G,\lambda'}(s,t)\le p_{H,\lambda}(s,t)$ and $c_{G,\lambda'}(s,t)\le c_{H,\lambda}(s,t)$. Indeed, these hold because the time-function $\lambda'$ does not allow for the existence of a \tpath{s,t} not contained in $H$.
	\end{proof}
	
	\section{Proof of sufficiency of Theorem~\ref{main}}\label{sec:sufficiency}

The following has been defined and proved in~\cite{KKK.00} and will be useful in our proof. Given a simple graph $G$, vertices $v,w\in V(G)$, and a positive integer $d$, a graph $G$ is called \emph{$(v,w,d)$-decomposable} if:

\begin{itemize}
\item Both $v$ and $w$ have degree $d$.
\item Either $G- \{v, w\}$ consists of $d$ components or $vw\in E(G)$ and $G-\{v,w\}$ has $d-1$ components.
\end{itemize}

\begin{lemma}[\cite{KKK.00}]\label{lem:KKK.00}
Let $G$ be a 2-connected simple graph with $\Delta(G)\ge 4$.  If $\mathcal{F}_{3}\not\preceq G$, then $G$ is $(v,w,d)$-decomposable for some $v, w \in V(G)$ and integer $d \geq 4$. 
\end{lemma}

\begin{figure}
    \centering
   \begin{subfigure}[b]{0.4\textwidth}
         \centering
        \includegraphics[width = 4.8cm]{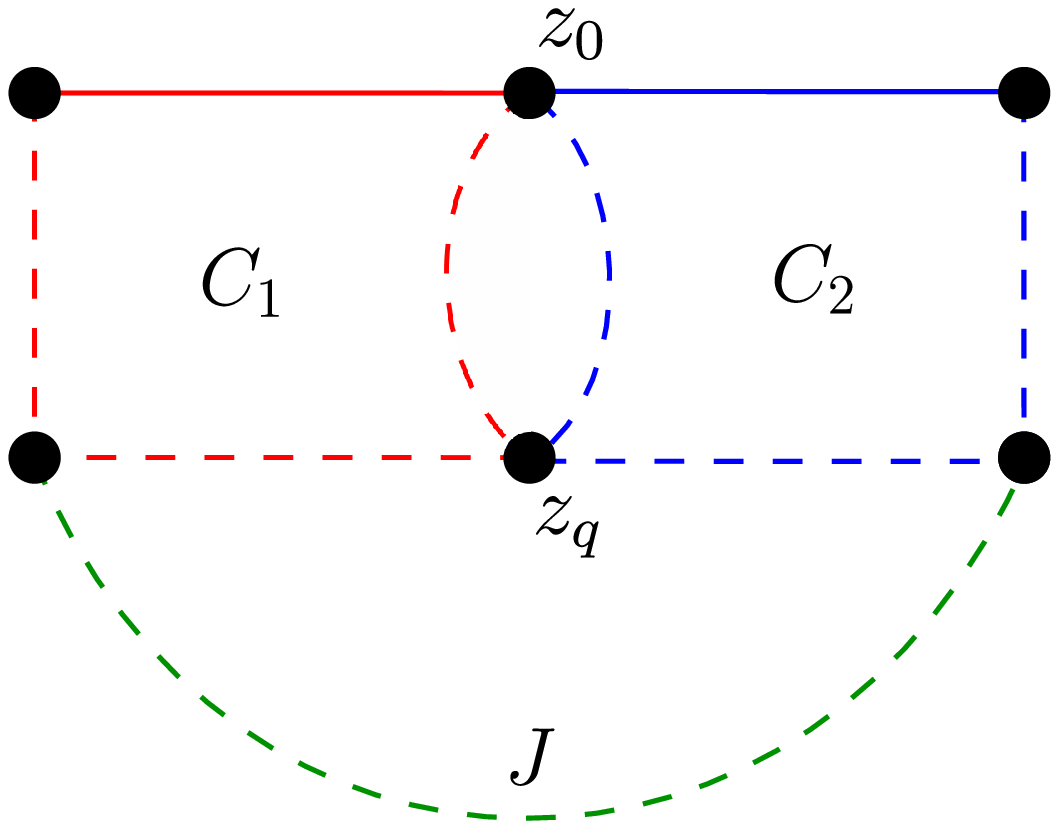}
        \caption{Cycle $C_1$ is highlighted in red, $C_2$ in blue, and path $J$ in green. \phantom{~~~~~~~~~~~~~~}}
        \label{fig:propf1}
     \end{subfigure}
   \hfill
   \begin{subfigure}[b]{0.5\textwidth}
        \centering
        \includegraphics[width = 4.8cm]{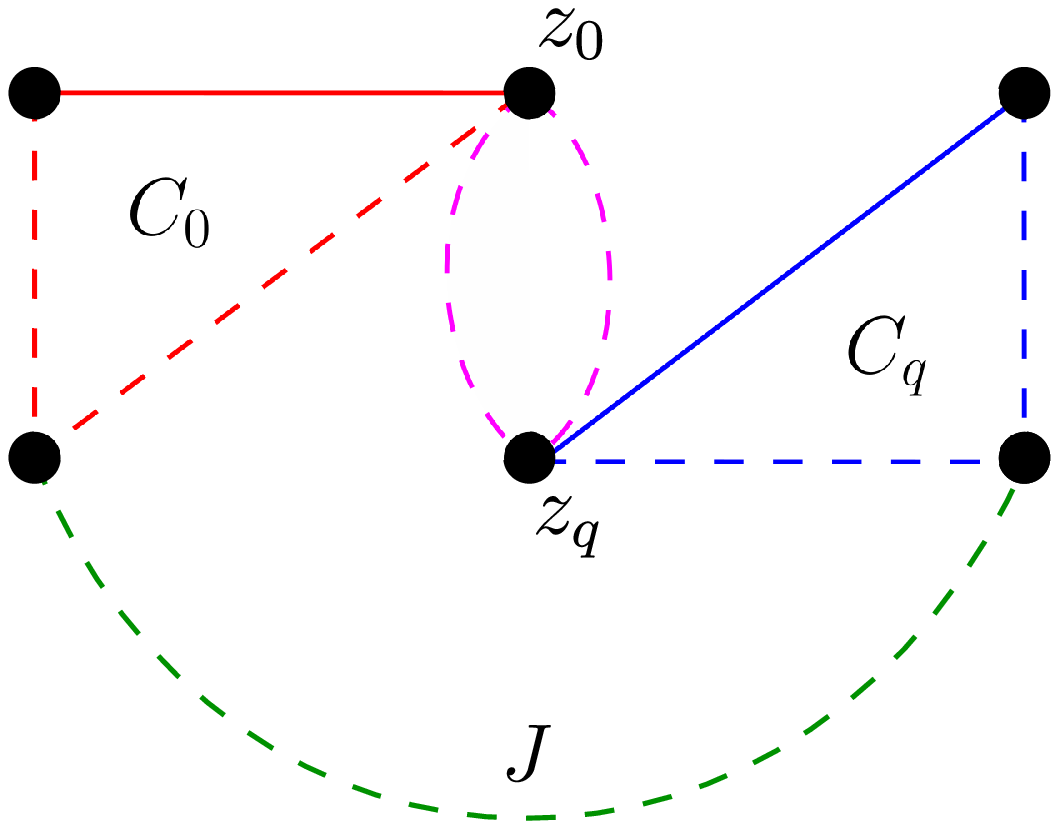}
        \caption{Cycle $C_0$ is highlighted in red, $C_q$ in blue, path $J$ in green, and the chain $L$ in magenta.}
        \label{fig:propf2}
     \end{subfigure}
     \caption{Dashed lines denote paths, dash parallel lines between $z_0$ and $z_q$ denote chains and solid lines denote  edges. }
\end{figure}

In our proof, we assume the existence of a counter-example and explore the properties of a minimum one in order to arrive to a contradiction. A \emph{minimum counter-example} is a graph $G$, together with $s,t\in V(G)$ and a time-function $\lambda$, that minimizes $\lvert V(G)\rvert + \lvert E(G)\rvert$ and is such that $p_{G,\lambda}(s,t) < c_{G,\lambda}(s,t)$ and $G$ has no $\mathcal{F}_i$ as an m-topological minor for every $i\in [3]$. 
In other words, $G$ is a minimum non-Mengerian graph that does not have $\mathcal{F}_i$ as an m-topological minor for every $i\in [3]$, while $\lambda,s,t$ are certificates that $G$ is non-Mengerian. 
The following straightforward propositions are useful to get to a contradiction by helping finding m-subdivisions of ${\cal F}_1$ and ${\cal F}_2$. 
Figures~\ref{fig:propf1} and~\ref{fig:propf2} help verification.

\begin{proposition}\label{prop:findF1}
Let $G$ be a graph, $L = (z_0,\cdots,z_q)$ be a chain in $G$, $C_1$ and $C_2$ be edge disjoint cycles and $J$ be a $C_1,C_2$-path. If all the items below hold, then the set of vertices and edges of $J, L, C_1, C_2$ form a subgraph of $G$ that is an m-subdivision of ${\cal F}_1$.

\begin{itemize}
    \item $V(C_1)\cap V(C_2) = V(L)$.
    \item $V(L)\cap V(J)=\emptyset$.
    \item Letting $w_1,w_2$ be the extremities of $J$, we have that $\{w_1,w_2\}\cap N(z_i)=\emptyset$ for some $i\in \{0,q\}$.
\end{itemize}
Reciprocally, if $H$ is an m-subdivision of $\mathcal{F}_1$ in $G$, then such $J,L,C_1,C_2$ exist that form $H$.
\end{proposition}

\begin{proposition}\label{prop:findF2} 
Let $G$ be a graph, $L = (z_0,\cdots,z_q)$ be a chain in $G$, $C_0$ and $C_q$ be vertex disjoint cycles in $G$ such that, for each $i\in \{0,q\}$, $V(C_i)\cap V(L) = \{z_i\}$, and $J$ be a $C_0,C_q$-path such that $V(L)\cap V(J)=\emptyset$. Then the set of vertices and edges of $L,J,C_0,C_q$ form a subgraph of $G$ that is an m-subdivision of ${\cal F}_2$. Reciprocally, if $H$ is an m-subdivision of $\mathcal{F}_2$ in $G$, then such $J,L,C_0,C_q$ exist that form $H$. 
\end{proposition}

\begin{figure}[!h]
   \centering
   \begin{subfigure}[b]{0.4\textwidth}
         \centering
         \includegraphics[width=4cm]{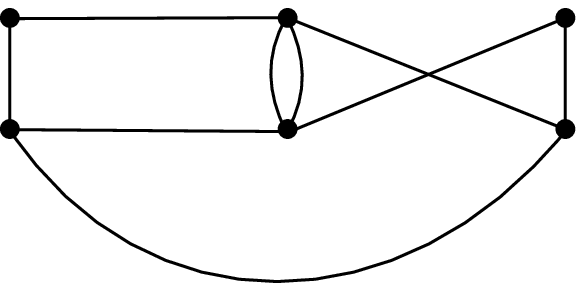}
         \caption{Smallest graph having a crossed structure.}
         \label{fig:crossed_graph}
     \end{subfigure}
   \hfill
   \begin{subfigure}[b]{0.5\textwidth}
         \centering
         \includegraphics[width=5.5cm]{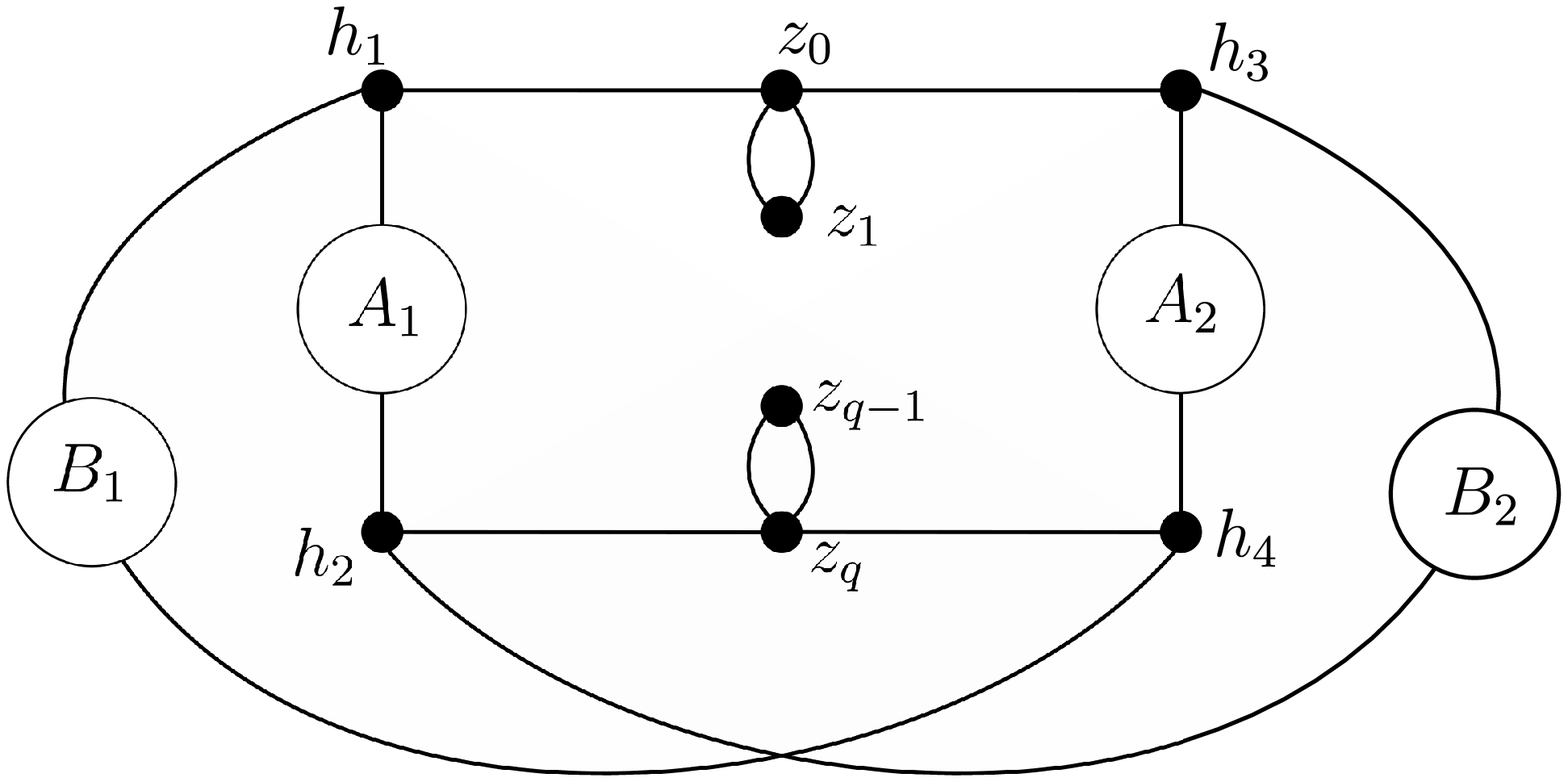}
         \caption{Crossed structure.}\label{fig:crossed_structure}
         \label{fig:five over x}
     \end{subfigure}
    \caption{Representation of a graph that has a crossed structure.}\label{fig:deg3}
\end{figure}

The general idea of the proof is to get to a contradiction by finding an m-subdivision of some $\mathcal{F}_i$ in a minimum counter-example. However, instead of the desired m-subdivision, sometimes we can only find an m-subdivision of the graph in Figure~\ref{fig:crossed_graph}, which we call \emph{crossed}. We say that $G$ has a \emph{crossed structure} if $V(G)$ can be partitioned as depicted in Figure~\ref{fig:crossed_structure}, where solid lines represent multiedges and each $X\in \{A_1,A_2,B_1,B_2\}$ represents either a multiedge  or a subgraph linked to the rest of the graph only through the depicted multiedges; see Figure~\ref{fig:crossed_graph} for the case where every $X\in\{A_1,A_2,B_1,B_2\}$ is a multiedge. Additionally, $B_2$ always exists (i.e., $h_2,h_3$ are always linked through a multiedge or through subgraph $B_2$), while $B_1$ might not exist (i.e., $h_1,h_4$ might not be linked through a multiedge or subgraph $B_2$). We call any $h_2,h_3$-path (or $h_1,h_4$-path) passing through $B_2$ ($B_1$) a \emph{crossing path}; also, we say that $G$ is \emph{1-crossed} if only $h_2,h_3$ have a crossing path, and \emph{2-crossed} otherwise. The next lemma is also useful for our recognition algorithm.

\begin{lemma}
Let $G$ be a 2-connected graph that has no $\mathcal{F}_3$ as topological minor, and let $L$ be a chain of $G$ for which all internal vertices have degree~2 in $U(G)$. If the graph obtained from $G$ by identifying the vertices of $L$ has a subdivision of $\mathcal{F}_3$, then one of the following holds, and we can decide which one in time $O(n^3)$, where $n = \lvert V(G)\rvert$:
\begin{itemize}
    \item $\mathcal{F}_1$ or $\mathcal{F}_2$ is an m-topological minor of $G$.
    \item $G$ has a crossed structure. 
\end{itemize}
\label{lem:helpcrossed}
\end{lemma}

\begin{proof}
Write $L$ as $(z_0,\ldots, z_q)$. First call the graph obtained from $G$ by identifying the vertices of $L$ by $G_L$ and let $H$ be a subdivision of $\mathcal{F}_3$ contained in $G_L$; denote by $\ell$ the vertex obtained by the identification of $L$. 
Such $H$ can be found in time $O(n^3)$ using the algorithm presented in~\cite{GKMW.11}. 
Note that $\ell \in V(H)$ as otherwise we would have that $H$ is a subgraph of $G$, a contradiction to the fact that $G$ has no $\mathcal{F}_3$ as topological minor. Note also that $N_{H}(\ell) = (N_G(z_0)\cup N_G(z_q))\setminus V(L)$ as all internal vertices of $L$ have degree~2 in $U(G)$, and that we can suppose $N_G(z_0)\cap N_G(z_q) = \emptyset$ as otherwise the removal of some edge incident to $z_0$ or $z_q$ would lead to the same graph $G_L$ and the argument below can still be applied. First we prove that $|N_{H}(\ell)\cap N_G(z_0)| = 2$. For this, we analyse the cases:

\begin{itemize}
    \item If $|N_{H}(\ell)\cap N_G(z_0)| = 0$, then $N_{H}(\ell)\subseteq N(z_q)$. Therefore the subgraph contained in $G$ obtained from $H$ by replacing $L$ by $z_q$ is isomorphic to $H$, a contradiction as $G$ has no subdivision of $\mathcal{F}_3$.
    
    \item If $|N_{H}(\ell)\cap N_G(z_0)| = 1$, then let $Q$ be a $z_0,z_q$-path in $G$ contained in $L$ and $a \in N_{H}(\ell)\cap N_G(z_0)$. Then, if we substitute the edge $a\ell$ in $H$ by the $a,z_q$-path $az_0Qz_q$, we would obtain a subgraph of $G$ that is a subdivision of $H$, again a contradiction. 
    
    \item If $|N_{H}(\ell)\cap N_G(z_0)| = 3$. Since $N_G(z_0)\cap N_G(z_q) = \emptyset$, $N_{H}(\ell) = (N_G(z_0)\cup N_G(z_q))\setminus V(L)$ and $\Delta(H)  =4$, we get that $|N_{H}(\ell)\cap N_G(z_q)| \le 1$, a contradiction as the previous cases can be analogously applied to $z_q$.
    
\end{itemize}

Therefore $|N_{H}(\ell)\cap N_G(z_0)| = 2$ and, by applying the same reasoning we get $|N_{H}(\ell)\cap N_G(z_q)| = 2$. 
Now let $H^*$ be the subgraph of $G$ formed by $L$ and $V(H)\cap V(G)$. Let $h_{1},\ldots,h_4$ be the vertices of $H^*$ corresponding to the vertices of degree~3 in  $\mathcal{F}_3$. Also, for each $i \in \{1,2,3,4\}$, let $P_i$ be the path in $G$ between $h_i$ and $N_{H}(\ell)$, and call by $u_i$ the extremity of $P_i$ distinct from $h_i$. 
We can find these vertices and paths by running a search on $H^*$, which takes $O(m+n)$ time where $m= \lvert E(G)\rvert$. 
As we have proved, in $G$ two of the vertices in $\{u_{1},\ldots,u_4\}$ are incident to $z_0$ and the other two are incident to $z_q$. All the possible adjacencies are depicted in Figure~\ref{fig:deg3p}.

\begin{figure}[ht]
   \centering
   \begin{subfigure}[b]{0.3\textwidth}
         \centering
         \includegraphics[width=3.6cm]{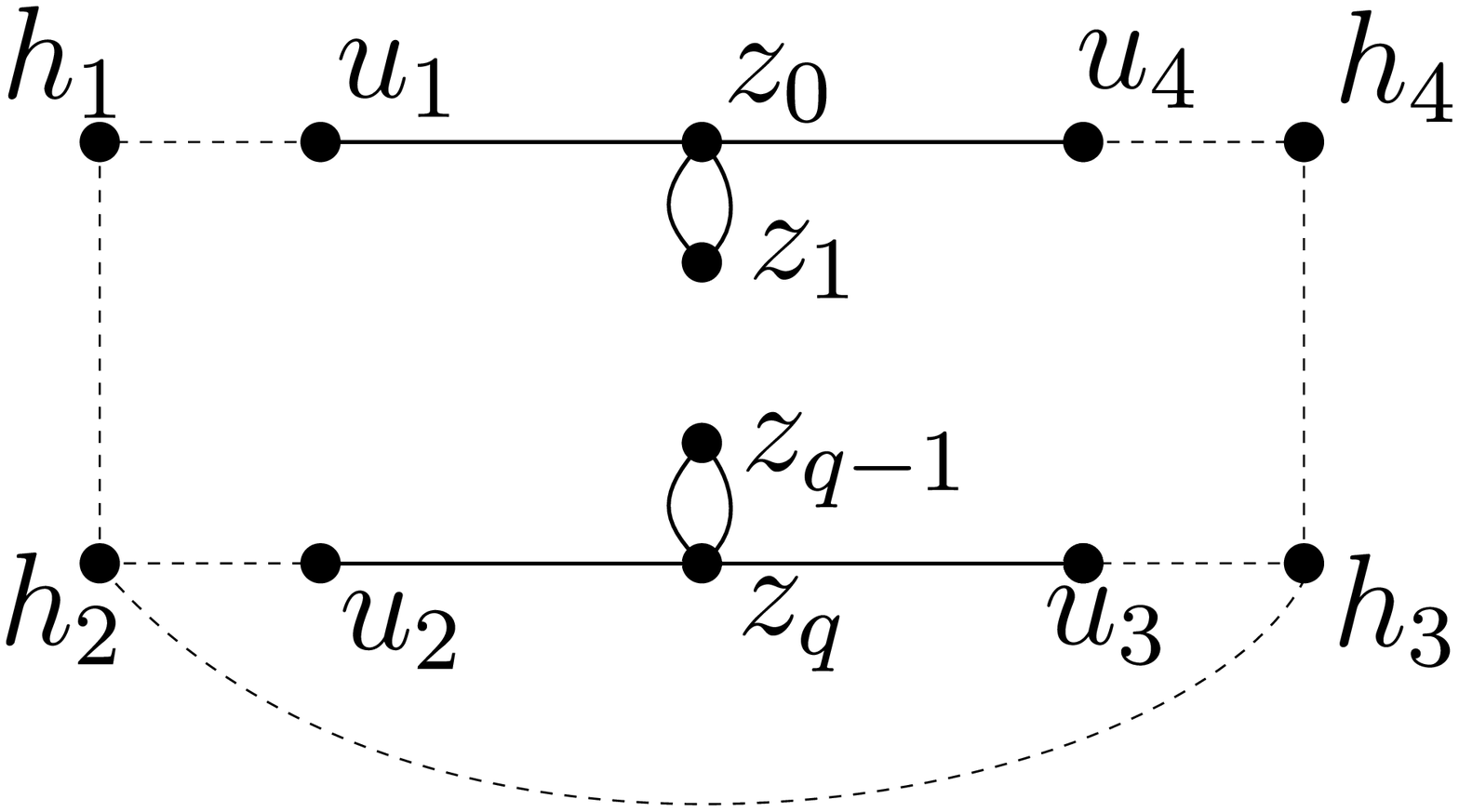}
         \caption{$u_1,u_2$ incident to $z_0$.}\label{fig:F1}
     \end{subfigure}
   \hfill
   \begin{subfigure}[b]{0.3\textwidth}
         \centering
         \includegraphics[width=3.6cm]{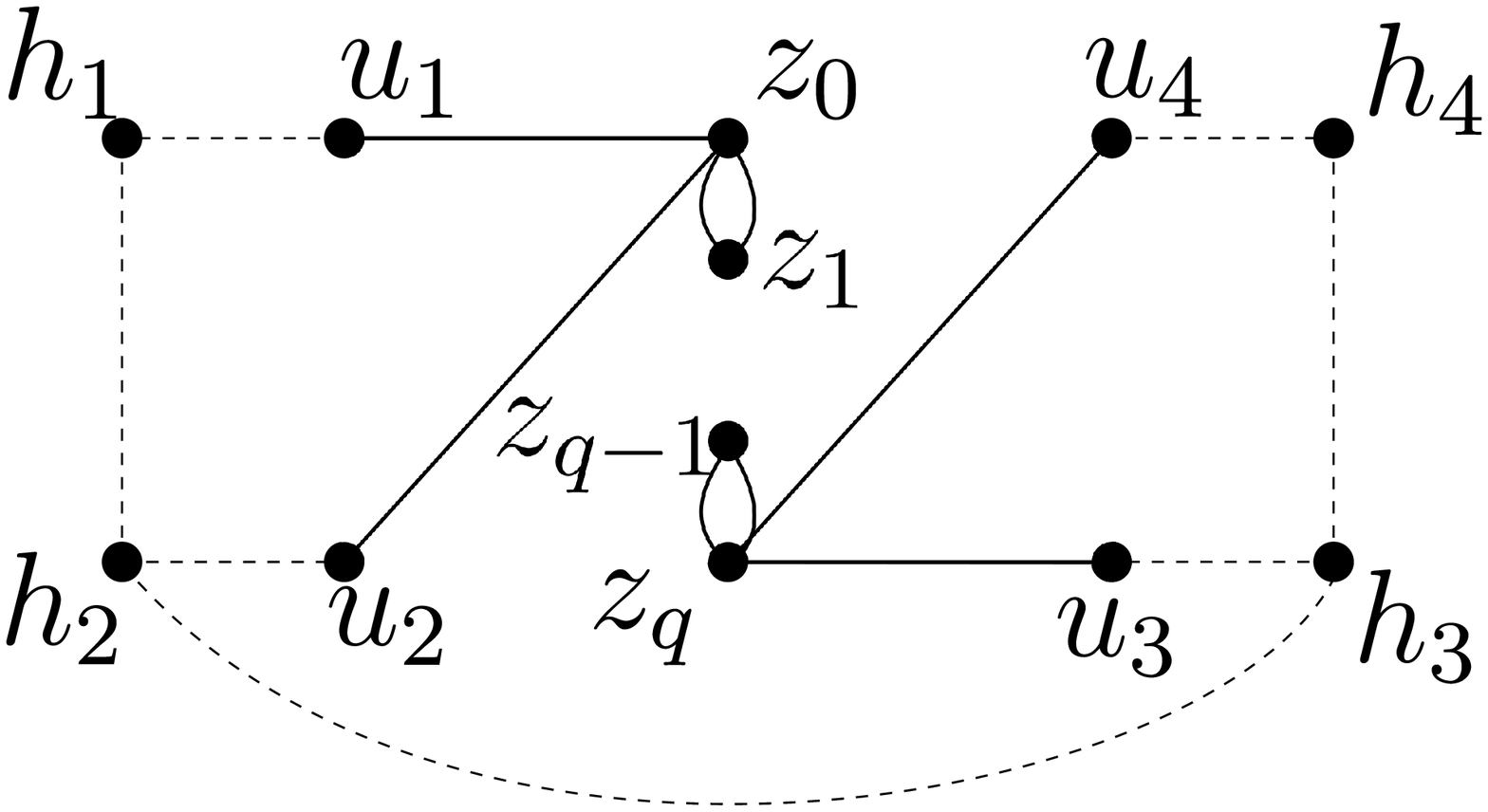}
         \caption{$u_1,u_4$ incident to $z_0$.}\label{fig:F3}
     \end{subfigure}
   \hfill
   \begin{subfigure}[b]{0.3\textwidth}
         \centering
         \includegraphics[width=3.6cm]{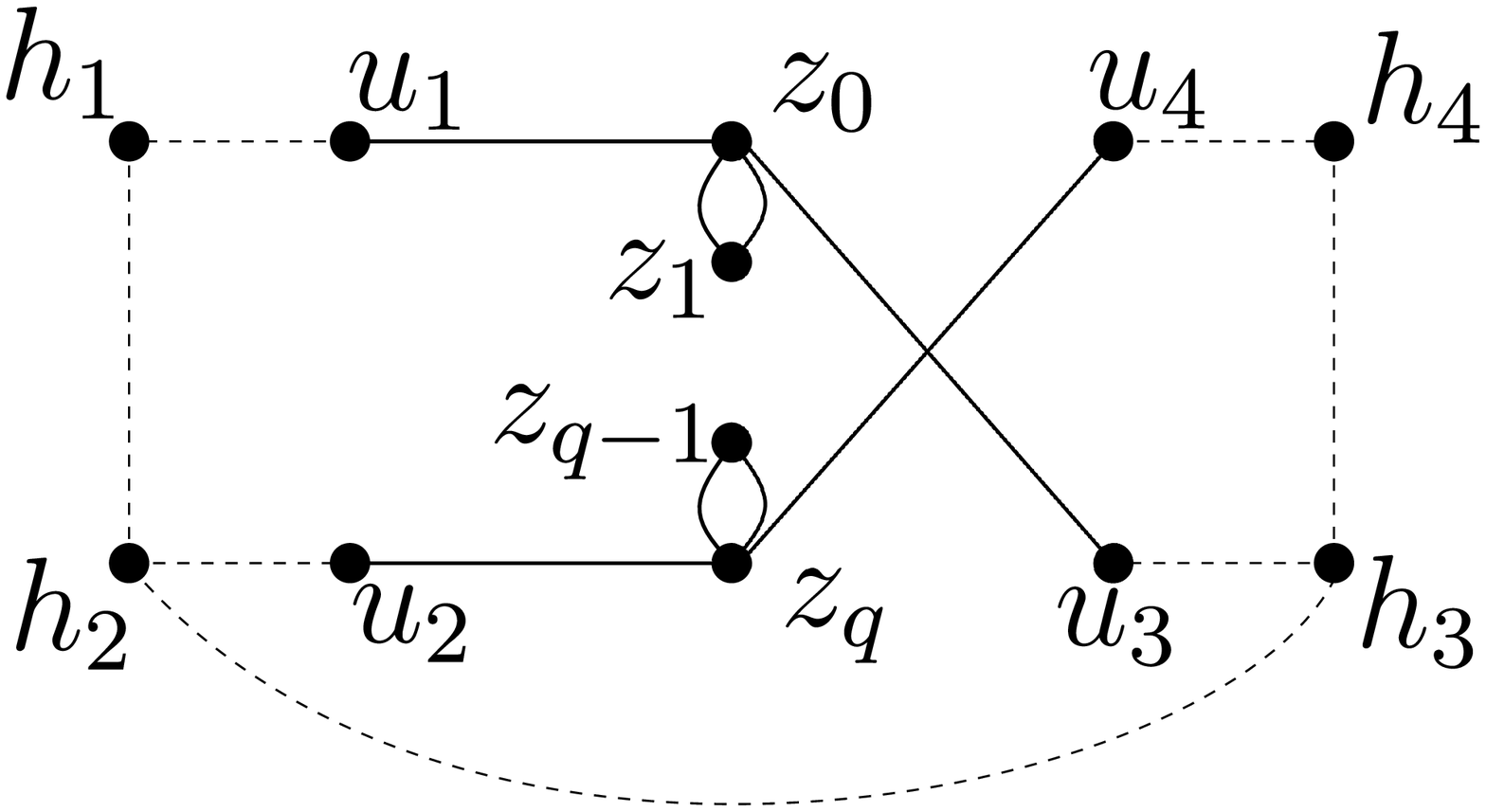}
         \caption{$u_1,u_3$ incident to $z_0$.}\label{fig:crossed_neighborhood}
     \end{subfigure}
    \caption{Possible adjacencies of $u_1,\cdots, u_4$ in $\{z_0,z_q\}$. Dashed lines denote paths, while solid lines denote  edges. The cases $u_1,u_2$ incident to $z_q$, $u_1,u_3$ incident to $z_q$, and $u_1,u_4$ incident to $z_q$ are  analogous.}\label{fig:deg3p}
\end{figure}

Observe that if the situation in Figure~\ref{fig:F1} occurs, then $H^*$ is an m-subdivision of ${\cal F}_1$, while if the situation in Figure~\ref{fig:F3} occurs, then $H^*$ is an m-subdivision of ${\cal F}_2$. 
Additionally, note that if $u_2 \neq h_2$ or $u_3 \neq h_3$, then $H^*$ is also an m-subdivision of $\mathcal{F}_1$. 
Since all these tests can be done in constant time, we can now consider that the situation depicted in Figure~\ref{fig:crossed_neighborhood} occurs and $h_2z_q$ and $h_3z_0$ are edges of $H^*$. 

Define $A_1$ to be the set of internal vertices in the $h_2,z_0$-path in $H^*$ passing by $h_1$; $A_2$ to be the set of internal vertices in the $h_3,z_q$-path in $H^*$ passing by $h_4$; and $B_1$ the set of internal vertices in the $h_2,h_3$-path not passing by $h_1$ or $h_4$. Also, for each $X \in \{A_1,A_2,B_1\}$, denote by $P_{X}$ the path formed by $X$. 
Notice that all vertices in $V(H^*)\setminus(A_1\cup A_2\cup B_1)$ either have degree~3 in $H^*$ or are internal vertices of the chain, which we know that have degree~2 in $U(G)$. Now we show that either both items below hold, or we can find an m-subdivision of $\mathcal{F}_1$. Note that this finishes our proof since this test can be done in time $O(m+n)$, and these items imply that $G$ has a crossed structure.

\begin{itemize}
    \item All paths between $B_{1}$ and $A_1\cup A_2$ in $G$ are contained in $H^*$.

    \item If there is a path between $A_1$ and $A_2$ in $G$ that is not contained in $H^*$, then it must be a path between $u_1$ and $u_4$. 
\end{itemize}

Suppose first that there exists a path $P'$ in $G$ between $a \in A_1$ and $b \in B_1$ not contained in $H^*$. We obtain the elements necessary to apply Proposition~\ref{prop:findF1}; observe Figure~\ref{fig:crossedlemma1} to follow the argument. Let $P'' = aP'bP_{B_1}h_3$. Also, let $C_1,C_2$ be two edge-disjoint cycles formed by $A_1\cup \{h_2\}$ and $A_2\cup \{h_3\}$ together with the chain $L$. 
We can apply Proposition~\ref{prop:findF1} to $L,P'',C_1,C_2$ to obtain an m-subdivision of ${\cal F}_1$. An analogous argument can be applied on a path between $B_1$ and $A_2$.

\begin{figure}[!h]
   \centering
   \begin{subfigure}[b]{0.4\textwidth}
         \centering
         \includegraphics[width=4.5cm]{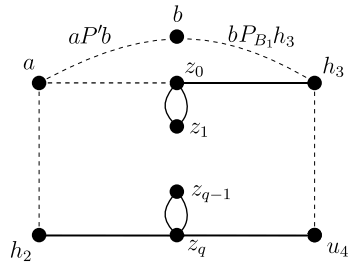}
         \caption{Subdivision of ${\cal F}_1$ when there exists $P'$ between $B_1$ and $A_1$ not in $H^*$.}
         \label{fig:crossedlemma1}
     \end{subfigure}
   \hfill
   \begin{subfigure}[b]{0.4\textwidth}
         \centering
         \includegraphics[width=4.5cm]{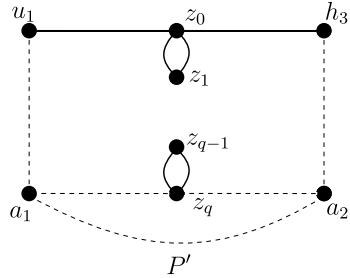}
         \caption{Subdivision of ${\cal F}_1$ when there exists $P'$ between $A_1\setminus\{u_1\}$ and $A_2$ not in $H^*$.}
         \label{fig:crossedlemma2}
     \end{subfigure}
    \caption{Case analysis in the proof that $G$ has a crossed structure. Dashed lines represent paths, while solid lines represent edges.}\label{fig:deg3q}
\end{figure}

Now, suppose that there is a path $P'$ between $a_1\in A_1$ and $a_2\in A_2$ in $G$ not contained in $H^*$, and suppose that $a_1\neq u_1$. Observe Figure~\ref{fig:crossedlemma2} to see that, again, we can apply Proposition~\ref{prop:findF1} to obtain an m-subdivision of $\mathcal{F}_1$. The same argument can be applied in case $a_2\neq u_4$ and therefore the second item above follows.

\end{proof}

In what follows, we aim to prove that either $\Delta(U(G))\ge 4$, or $G$ has a crossed structure. For this, the following lemma will be useful. 

\begin{lemma}
\label{lemma:simple_st}
Let $(G,\lambda)$, $s,t$ be a minimum counter-example.
Then, $G$ is 2-connected and every multiedge incident in $s$ or $t$ has multiplicity~1. 
\end{lemma}
\begin{proof}
	We first observe that $G$ is connected. If $s$ and $t$ are in distinct components, then $p_{G,\lambda}(s,t) = c_{G,\lambda}(s,t)= 0$. Then $s,t$ are in the same component and as $(G,\lambda)$ is minimum, it follows that $G$ is connected.  Now, we prove that $G$ is 2-connected. By contradiction, suppose that $\{v\}$ is a \blue{cut} of $G$, and $G_{1},\cdots, G_{d}$ be the subgraphs induced by the vertex set of each component of $U(G)-v$ and $v$. If $s$ and $t$ belong to $G_{i}$ for some $i \in [d]$, then every \tpath{s,t} is contained in $G_{i}$, making $G_{i}$ a non-Mengerian proper subgraph of $G$, contradicting the minimality of $G$. So suppose that $s \in G_{i}$ and $t \in G_{j}$ for $i,j \in [d]$, with $i \neq j$. Observe that in this case, either there is no \tpath{s,t}, in which case the empty set is a \vcut{s,t} (i.e. $p_{G,\lambda}(s,t) = c_{G,\lambda}(s,t) = 0$), or $p_{G,\lambda}(s,t) = c_{G,\lambda}(s,t) = 1$ as every \tpath{s,t} goes through $v$. This contradicts the definition of minimum counter-example.

	Now, by contradiction, suppose that $sy\in E(U(G))$ has multiplicity at least~$2$ in $G$. Let $e\in E(G)$ be the edge with endpoints $sy$ that minimizes $\lambda(e)$, and let $S = \{e'\in E(G)\mid e'\neq e \mbox{ and } e'\mbox{ has endpoints }sy\}$. Denote by $(G^*,\lambda)$ the temporal graph $(G-S,\lambda)$. Because $(G^*,\lambda)$ is contained in $(G,\lambda)$, 
	observe that we get that $p_{G^*,\lambda}(s,t)\le p_{G,\lambda}(s,t)$ and $c_{G^*,\lambda}(s,t)\le c_{G,\lambda}(s,t)$. 
	In addition, notice that a set ${\cal P}$ of vertex disjoint temporal paths in $(G,\lambda)$ also corresponds to a set of vertex disjoint temporal paths in $(G^*,\lambda)$ of same size, since it suffices to replace, if needed, the first edge $e'\in S$ with endpoints $sy$ by $e$. This gives us that $p_{G,\lambda}(s,t) = p_{G^*, \lambda}(s, t)$ and $c_{G, \lambda}(s, t) = c_{G^*, \lambda}(s, t)$, contradicting the minimality of $G$. A similar argument can be applied if some edge incident in $t$ has multiplicity at least~2, with the difference that we take its  latest occurrence instead. 
\end{proof}

\begin{lemma}
Let $(G,\lambda)$, $s,t$ be a minimum counter-example. If $\Delta(U(G)) \leq 3$, then $G$ has a crossed structure. 
\label{lem:crossed}
\end{lemma}
\begin{proof}
We start by finding the chain of $G$.
 Denote by $p'_{G,\lambda}(s,t)$ the maximum number of edge disjoint \tpath{s,t}s, i.e., the maximum $k$ such that there are \tpath{s,t}s $P_1,\ldots,P_k$ with $E(P_i)\cap E(P_j)=\emptyset$ for every $i,j\in [k]$, $i\neq j$. Define also $c'_{G,\lambda}(s,t)$ as the minimum cardinality of a subset $S\subseteq E(G)$ such that $(G-S,\lambda)$ has no \tpath{s,t}s. 
 Note that, since $st\notin E(U(G))$, then removing one endpoint of each $e\in S$ also destroys every \tpath{s,t}; hence $c'_{G,\lambda}(s,t) \geq c_{G,\lambda}(s,t) > p_{G,\lambda}(s,t)$. Using the equality $c'_{G,\lambda}(s,t) = p'_{G,\lambda}(s,t)$ proved in~\cite{B.96}, one can find at least $c_{G,\lambda}(s,t)$ edge disjoint \tpath{s,t}s. Because $c_{G,\lambda}(s,t) > p_{G,\lambda}(s,t)$, at least two of these paths, say $J_1,J_2$, must intersect in an internal vertex. 
 Now, let $z_0$ be the closest to $s$ in $(V(J_1)\cap V(J_2))\setminus \{s,t\}$. For each $i\in \{1,2\}$, let $y_i$ be the neighbor of $z_0$ in $sJ_iz_0$. Since $d_{U(G)}(z_0)\le 3$, $z_0\neq t$ and $y_1\neq y_2$ by the choice of $z_0$, it follows that $d_{U(G)}(z_0)=3$ and the next vertex in $J_1$ and $J_2$ also coincides, i.e., that there exists $z_1\in (N(z_0)\cap V(J_1)\cap V(J_2))\setminus \{y_1,y_2\}$. Now consider $(z_0,z_1,\cdots, z_q)$ be maximal such that $(z_0,e_1,z_1,\cdots,e_q,z_q)$ is a subpath in $J_1$ and $(z_0,e'_1,z_1,\cdots,e'_q,z_q)$ is a subpath in $J_2$. Since $J_1$ and $J_2$ are edge disjoint, we get that $e_i\neq e'_i$  for every $i\in [q]$. Additionaly, by Lemma~\ref{lemma:simple_st}, we have $z_q\neq t$, and therefore $d_{U(G)}(z_q) = 3$. 
 
Now let $L = (v = z_{0},\cdots, z_{q} = w)$ be the chain previously obtained, and suppose that $d_{U(G)}(z_{i}) = 2$ for every $i \in [q-1]$; otherwise, it suffices to shorten the chain until it is appropriate. Let $G'$ be obtained by the identification of $L$, and let $\ell$ denote the new vertex. Notice that $U(G')$ is $2$-connected and has exactly one vertex of degree at least $4$, namely $\ell$.  Observe that Lemma~\ref{lem:KKK.00} tells us that every $2$-connected graph with exactly one vertex of degree at least~4 must have a subdivision of ${\cal F}_3$. Recall that $G$ does not have $\mathcal{F}_i$ as m-topological minor for each $i\in [3]$, as $(G,\lambda),s,t$ is a counter-example. It follows from Lemma~\ref{lem:helpcrossed} that $G$ has a crossed structure.
\end{proof}

In order to treat the case where $\Delta(U(G))\le 3$, we prove that a graph $G$ that has a crossed structure is in fact Mengerian, and hence cannot compose a minimum counter-example. For this, the following structural theorem will be useful.

\begin{lemma}\label{lem:edgecut}
Let $(G,\lambda)$, $s,t$ be a minimum counter-example. Suppose that $S = \{uu',vv'\}\subseteq E(U(G))$ is such that $U(G) - S$ is disconnected, and consider that $u',v'$ are in the same component $G'$ of $U(G)-S$. If $\min\{d_{U(G)}(u),d_{U(G)}(v)\} \geq 3$, then:
    
    \begin{enumerate}
        \item\label{lem:edgecut.i} At least one of $s,t$ is not in $G'$.
        
        \item\label{lem:edgecut.ii} If $x \in V(G')\cap \{s,t\}$, then $V(G') = \{x\}$. 
        
        \item\label{lem:edgecut.iii} If $s,t \notin V(G')$, then $G'$ is an $u',v'$-path.
    \end{enumerate}
\end{lemma}
\begin{proof}
In what follows, for simplicity, we denote $d_{U(G)}(x)$ by $d(x)$. We divide the proof into cases.

\paragraph{Case 1.} $uu'$ and $vv'$ have multiplicity $1$. 
If $\{s,t\} \subseteq V(G')$, then let $(G^*,\lambda^*)$ be obtained by identifying all vertices in $V(G)\setminus V(G')$. Observe that $\lvert V(G^*)\rvert < \lvert V(G)\rvert$ (as $V(G)\setminus V(G')$ contains at least two vertices, $u$ and $v$). Additionally, since $uu'$ and $vv'$ have multiplicity~1 and every \tpath{s,t} not contained in $G'$ must use both $uu'$ and $vv'$, observe that this operation preserves the \tpath{s,t}s not contained in $(G',\lambda)$. Finally, note that this also cannot create new m-topological minors as the vertex obtained by the contraction has degree 2 in $G^*$. Therefore, we get that $(G^*,\lambda^*),s,t$ is a smaller counter-example, a contradiction. 
Observe that, in case $\{s,t\} \cap V(G')=\emptyset$, then either $u'=v'$ and Item~\ref{lem:edgecut.iii} follows, or a similar argument can be applied  by identifying all vertices of $V(G')$. We  therefore get that $\lvert\{s,t\}\cap V(G')\rvert=1$.

Now, suppose that $s\in V(G')$ and $t \notin V(G')$. In what follows, we apply again identification operations that produce vertices of degree 2; hence no m-topological minors are created, and we consider this implicitly. Observe that, as $\{u,v\}$ is a \vcut{s,t}, we have $1 = p_{G,\lambda}(s,t) < c_{G,\lambda}(s,t) = 2$. So call the only edge with endpoints $uu'$ by $e$ and the only edge with endpoints  $vv'$ by $f$. 
First, suppose that (*) every two \tpath{s,t}s $Q_{1}, Q_2$ intersect in $V(G')\setminus \{s\}$, and  
let $(G^*,\lambda^*)$ be the temporal graph obtained by identifying all vertices in $Z = V(G)\setminus V(G')$ into a single vertex $z$.
Observe that by (*) we get $p_{G^*,\lambda^*}(s,z) = 1$, and since $G^*$ is smaller than $G$, we then must have $c_{G^*,\lambda^*}(s,z)=1$. This is a contradiction as any vertex separating $s$ from $z$ in $(G^*,\lambda^*)$ would also separate $s$ from $t$ in $(G,\lambda)$. Thus there are two \tpath{s,t}s in $(G,\lambda)$, $Q_{1}$ and $Q_{2}$, such that $V(Q_1)\cap V(Q_2) \cap V(G') = \{s\}$. 
Finally, suppose that $\lvert V(G')\rvert >1$, and now let $(G^*,\lambda^*)$ be obtained from $(G,\lambda)$ by identifying all vertices in $V(G')$ into a single vertex $x$. From the existence of $Q_1,Q_2$ and  the fact that any \tpath{s,t} must use $e$ or $f$, observe that we must have $p_{G^*,\lambda^*}(x,t)=1$. Again because $G^*$ is smaller than $G$, we get $c_{G^*,\lambda^*}(x,t)=1$, a contradiction as any such vertex would also separate $s$ from $t$ in $(G,\lambda)$.

\paragraph{Case 2.} $uu'$ has multiplicity at least $2$. 
As $G$ is 2-connected (Lemma~\ref{lemma:simple_st}) and $d(u)\geq 3$, there is a path $P$ between $x,y\in N(u)\setminus \{u'\}$ that does no pass by $u$; notice that $P$ does not intersect $V(G')$ as $u,v$ separates $G'$ from $x,y$. Therefore, the cycle $C_u$ formed by $P$ and $u$ does not intersect $V(G')$. Now, if $d(u')\geq 3$, then we can proceed in the same way to obtain a cycle $C_{u'}$ contained in $G'$. Consider then a path $J_1$ in $G-u$ from $C_u-u$ to $v'$, and a path $J_2$  in $G'-u'$ from $v'$ to $C_{u'}-u'$, and observe that we can apply Proposition~\ref{prop:findF2} to get an m-subdivision of $\mathcal{F}_3$, a contradiction. Therefore $d(u') = 2$. Call $u_{2}$ the vertex in $N(u')\setminus \{u\}$. 
Observe that a similar argument can be applied to $u'u_2$ in case it has multiplicity at least~$2$, as long as we do not arrive to $v$. In other words, by letting $(u=u_0,u'=u_1,\ldots,u_q)$ be maximal such that $d(u_i)=2$ for every $i\in [q-1]$, we get that: 

\begin{enumerate}[label=(\alph*)]

    \item\label{lem:edgecut.case2.I} Either $(u_0,u_{1},\ldots,u_{q-1})$ is a chain and $u_{q-1}u_q$ has multiplicity~1; 

    \item\label{lem:edgecut.case2.II} Or $(u_0,u_{1},\ldots,u_{q})$ is a chain and $u_q = v$.
   
\end{enumerate}

Observe that if $vv'$ also has multiplicity at least~2, then we can also obtain a maximal chain $(v=v_0,v'=v_1,\ldots,v_p)$ such that $d(v_i)=2$ for every $i\in [p]$.
 
We now prove that  $s\notin V(G')$. By contradiction, suppose otherwise, and note that Lemma~\ref{lemma:simple_st} gives us Item~\ref{lem:edgecut.case2.I} above. Let $\alpha = \lambda(u_{q-1}u_q)$, and denote by $E'$ the set of edges with endpoints $u_{q-2}u_{q-1}$; we know that there are at least two such edges. Also observe that $\lambda(e)\neq \lambda(e')$ for every $e,e'\in E'$ with $e\neq e'$ as we are in a minimum counter-example. Over all $e\in E'$, let $e_1$ be the edge minimizing $\lambda(e)$, and $e_2$ the one maximizing $\lambda(e)$. If $\lambda(e_1) \ge \alpha$, then $\lambda(e_2) > \alpha$ and hence all the paths passing by $e_2$ must traverse $u_qu_{q-1}$ from $u_q$ to $u_{q-1}$ before using $e_2$; therefore, every such path could use the edge $e_1$ instead. In other words, $(G-\{e_2\},\lambda),s,t$ is still a counter-example, contradicting the minimality of $(G,\lambda)$.  We then get $\lambda(e_1) < \alpha$. 
A similar argument can be applied to conclude that $\alpha < \lambda(e_2)$. 

Now, note that every \tpath{s,t} using $e_1$ must use it right before $u_{q-1}u_q$, while one that uses $e_2$ must do so right after using $u_{q-1}u_q$. 
Because of this, we get that if $t\notin V(G')$, then a \tpath{s,t} $P$ passing by $e_1$ 
must first use some edge with endpoints $vv'$, then use some path contained in the chain $(u,u_1,\cdots,u_{q-1})$, a contradiction as in this case the other endpoint of $P$ cannot be $t$. 
Therefore, we must have $t\in V(G')$. 
Let $J_{1}$ be a \tpath{s,t} using $e_1$ and $J_{2}$ a \tpath{s,t} using $e_2$ (observe Figure~\ref{fig:pathsJ1J2}). Also, let $E_{1}$ be the set of edges in $J_{1}$ appearing before time $\alpha$ and $E_{2}$ be the set of edges in $J_{2}$ appearing after time $\alpha$. Because $\max \lambda(E_{1}) = \lambda(e_1) < \alpha < \min \lambda(E_2) = \lambda(e_2)$, these sets are disjoint. 
Now observe that $u,u_1,\cdots,u_q$ are all contained in both $J_1$ and $J_2$, and this implies that $J_1,J_2$ must also contain $v$ and $v'$. Note additionally that the $s,v$-path contained in $J_1$ uses only edges in $E_1$, while the $v,t$-path contained in $J_2$ uses only edges in $E_2$. Since $E_1,E_2$ are disjoint, observe that, by picking a maximal chain $(v=v_0,v'=v_1,\ldots,v_p)$ such that $\{v_0,v_1,\ldots,v_p\}\subseteq V(J_1)\cap V(J_2)$, we get that $v_p$ must have degree at least~3, contradicting what is said in the paragraph after Item~\ref{lem:edgecut.case2.II}. 
We then conclude that $s\notin V(G')$ and an analogous approach can be used to conclude that $t\notin V(G')$; observe that this gives us Items~\ref{lem:edgecut.i} and~\ref{lem:edgecut.ii}, and hence it just remains to prove Item~\ref{lem:edgecut.iii}.

\begin{figure}[!h]
   \centering
    \includegraphics[width=8cm]{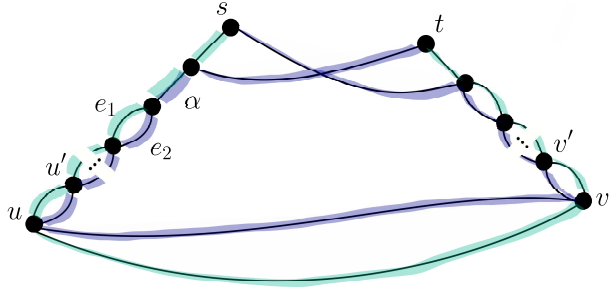}
    \caption{Paths $J_1$ and $J_2$. Subset $E_1$ is shadowed in blue, and $E_2$ is shadowed in green. Edge labeled with $\alpha$ has multiplicity~1 and is contained in both paths.}\label{fig:pathsJ1J2}
\end{figure}

Suppose by contradiction that Item~\ref{lem:edgecut.iii} does not hold, i.e., that $G'$ is not a $u',v'$-path. 
Observe that this implies that Item~\ref{lem:edgecut.case2.I} occurs as Item~\ref{lem:edgecut.case2.II} would imply that $G'$ is a $u',v'$-path. 
Now, we can apply the same argument to find a path $P = (v = v_1, v'=v_2,\ldots,v_{p-1},v_p)$ such that $(v_1,\ldots,v_{p-1})$ is either a chain or a single vertex, $d(v_i)=2$ for every $i\in [p-1]$, and $v_{p-1}v_p$ has multiplicity~1 (note that $v_{p-1}v_p$ can be equal to $vv'$). Let $e$ denote $u_{q-1}u_q$ and $f$ denote $v_{p-1}v_p$. 
Define  $X \colon = V(G')\setminus (\{u_1,\ldots,u_q\} \cup V(P))$; as $G'$ is not a $u'v'$-path, we get $X\neq \emptyset$. Let $(G^*,\lambda^*)$ be the temporal graph obtained by identifying $X\cup \{u_q,v_p\}$ into a single vertex, $z$. As $z$ is incident only to the edges corresponding to $e$ and $f$ in $G^*$, we get that this cannot have created an m-topological minor. Finally, observe that because $e$ and $f$ have multiplicity~1, we get that $p_{G^*,\lambda^*}(s,t) = p_{G,\lambda}(s,t)$ and $c_{G^*,\lambda^*}(s,t) = c_{G,\lambda}(s,t)$, a contradiction as $(G,\lambda),s,t$ is a minimum counter-example.
\end{proof}

The next lemma will allow us to avoid explicitly analyzing some cases in the final steps of our proof. Let $(G,\lambda)$ be a temporal graph and $T$ be its lifetime, i.e., the maximum value $\lambda(e)$ over all $e\in E(G)$. The \emph{reverse} of $(G,\lambda)$ is the temporal graph $(G,\lambda^-)$ where $\lambda^- \colon E(G) \to \mathbb{N} - \{0\}$ is defined below:

\[\lambda^-(e) = T + 1 - \lambda(e) .\]

\begin{lemma}\label{lem:reverse}
	Let $(G,\lambda)$ be a temporal graph and $s,t \in V(G)$. Then,
	
	\[c_{G,\lambda}(s,t) = c_{G,\lambda^-}(t,s) \text{ and } p_{G,\lambda}(s,t) = p_{G,\lambda^-}(t,s).\]
\end{lemma}

\begin{proof}
	It suffices to see that if $J = (s = v_1,e_1,\ldots,e_k,v_{k+1} = t)$ is a \tpath{s,t} in $(G,\lambda)$, then $J^- = (t = v_{k+1},e_{k+1},\ldots ,e_1,v_1 = s)$ is a \tpath{t,s} in $(G,\lambda^-)$, and vice-versa. As $V(J)= V(J^-)$, the result holds.
\end{proof}

Recall that, in order to conclude the case $\Delta(U(G))\le 3$, we want to prove that if $G$ has a crossed structure, then $G$ is Mengerian. This together with Lemma~\ref{lem:crossed} gives us that, in fact, a minimum counter-example cannot satisfy $\Delta(U(G))\le 3$. This is the idea of the proof of our final tool lemma.

\begin{lemma}
If $(G,\lambda),s,t$ is a minimum counter-example, then $\Delta(U(G))\geq 4$.\label{lem:Delta_geq_4}
\end{lemma}
\begin{proof}

By contradiction, suppose that $\Delta(G)\le 3$. By Lemma~\ref{lem:crossed}, we get that $G$ has a crossed structure. 
We use the same terminology as before for the vertices of $G$ (see~Figure~\ref{fig:crossed_structure}), i.e., $(z_0,\ldots,z_q)$ is the chain, $h_1,h_3\in N(z_0)$, $h_2,h_4\in N(z_q)$, $h_1,h_2$ are linked through $A_1$, $h_3,h_4$ are linked through $A_2$, $h_1,h_4$ are linked through $B_1$, and $h_2,h_3$ are linked through $B_2$, in case $G$ is 2-crossed. Additionally, when referring to the degree of a vertex, we always consider the degree in $U(G)$. 

We want to arrive to a contradiction by showing that $p(s,t) = c(s,t)$. Observe that if $p(s,t) = 0$, then it follows directly. So suppose that $p(s,t) \ge 1$ (and hence $c(s,t)\ge 2$). 
We first prove some further structural properties. 

\begin{claim}\label{claim:morestrucutre}The following hold.
\begin{enumerate}
    \item\label{claim:morestructure_item2} For every $x\in N(s)\cup N(t)$, we have that $d_{U(G)}(x) = 3$;
    
    \item\label{claim:morestructure_item1} For each edge $f \in E(G)$, there is a path $P$ using $f$ and not containing some vertex of degree~3 in $U(G)$; 

    \item\label{claim:morestructure_item3} Multiedges $z_qh_2$ and $z_0h_3$ have multiplicity~1. Additionally, if $G$ is 2-crossed, then also $z_0h_1$ and $z_qh_4$ have multiplicity~1.

    \item\label{claim:morestructure_item4} If $c(s,t) = 2$ and there is a vertex $v \in N(s)\cap N(t)$, then the multiedge incident to $v$ and $w\notin \{s,t\}$ have multiplicity at least~2.
\end{enumerate}
\end{claim}
\begin{proof}

For Item~\ref{claim:morestructure_item2}, suppose that $s$ has a neighbor $x$ such that $N(x) = \{s,y\}$ (i.e., $d(x) = 2$). Then every $e\in E(G)$ with endpoints $xy$ is such that $\lambda(e) \ge \lambda(sx)$, and we can remove $x$ and make $s$ directly adjacent to $y$ through an edge active at time $\lambda(sx)$ to obtain a smaller counter-example, a contradiction. Item~\ref{claim:morestructure_item2} follows by applying a similar argument to $t$. 

Now for Item~\ref{claim:morestructure_item1}, suppose by contradiction that $f\in E(G)$ is such that every \tpath{s,t} using $f$ contains all vertices of degree~3. Observe that the previous item gives us that every such path contains all neighbors of $s$. 
In other words, every path using $f$ contains a \vcut{s,t}, which in turn implies that the maximum number of disjoint \tpath{s,t}s in $(G-f,\lambda)$ is the same as in $(G,\lambda)$. Since $(G-f,\lambda)$ cannot be a counter-example, we get that there exists a \vcut{s,t} $S$ of $(G-f,\lambda)$ with less than $c_{G,\lambda}(s,t)$ vertices. 
If $S$ contains some vertex of degree~3, then $S$ also intersects every path using $f$, a contradiction as $S$ would be a \vcut{s,t} in $(G,\lambda)$. 
Hence consider $u \in S$, which must have degree~2. If $u$ lies in a \tpath{s,t} $P$ whose internal vertices all have degree~2, then by removing $P$ from $(G,\lambda)$ we decrease $p_{G,\lambda}(s,t)$ and $c_{G,\lambda}(s,t)$ by exactly one, a contradiction with the minimality of $(G,\lambda),s,t$. Therefore, every such path containing $u$ must contain some vertex of degree~3 different from $s$ and $t$; let $v$ be such a vertex. 
Note that in fact every \tpath{s,t} containing $u$ must also contain $v$. 
This means that $T = (S\setminus \{u\}) \cup \{v\}$ is also a \vcut{s,t} in $(G -f ,\lambda)$. Moreover, as every \tpath{s,t} that uses $f$ also must contain every vertex of degree~3, and hence in particular $v$, we get that $T$ is a \vcut{s,t} in $(G,\lambda)$, a contradiction as $|T| = |S| < c_{G,\lambda}(s,t)$.

For Item~\ref{claim:morestructure_item3}, recall Proposition~\ref{prop:findF1} and observe Figure~\ref{fig:h2h3multiple} to see that if $z_qh_2$ have multiplicity at least~2, then we can find an m-subdivision of ${\cal F}_1$; note also that the analogous holds for $z_0h_3$, and for $z_0h_1$ and $z_qh_4$ as well in case $G$ is 2-crossed. 

\begin{figure}[thb]
\begin{center}
\includegraphics[width=6cm]{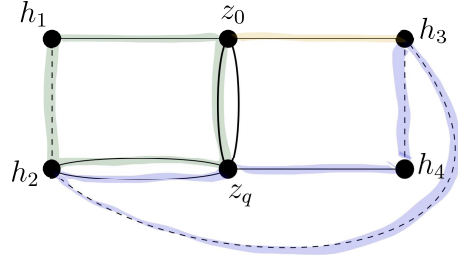}
\caption{Structure of $G$ in case $h_2z_q$ has multiplicity at least~2. The chain is simplified to two parallel edges. Dashed lines denote paths and solid lines denote edges.}
\label{fig:h2h3multiple}
\end{center}
\end{figure}

Finally, for Item~\ref{claim:morestructure_item4} suppose by contradiction that $vw$ has multiplicity~1 (note also that $w$ is well defined because of Item~\ref{claim:morestructure_item2} and the fact that $\Delta(U(G)) = 3$). Let $\alpha = \lambda(sv)$, $\beta = \lambda(vt)$ and $\gamma = \lambda(vw)$. 
If $\beta < \alpha$, then every path using $sv$ must use $vw$ and we have $\alpha \leq \gamma$; in the same way every path using $vt$ must use $vw$ and we have $\gamma \leq \beta$, a contradiction. Therefore we must have $\alpha \leq \beta$. Additionally, observe that any \tpath{s,t} in $(G - v,\lambda)$ is vertex disjoint from $(s,\alpha,v,\beta,t)$, a contradiction as in this case either $p(s,t)$ is also equal to~2, or $\{v\}$ is \vcut{s,t}. 
\end{proof}

We are now ready to prove the lemma. First note that Lemma~\ref{lemma:simple_st} gives us that $\{s,t\}\cap \{z_0,\ldots,z_q\} = \emptyset$. We now analyse the cases. 
    	
\begin{itemize}
	\item $G$ is 1-crossed. Note that Lemma~\ref{lem:edgecut} gives us that $A_1$ ($A_2$) is either the multiedge $h_1h_2$ ($h_3h_4$) or is a $h_1,h_2$-path ($h_3,h_4$-path). Hence both $h_2$ and $h_3$ have neighbors whose degree equals~2 in $U(G)$, and by Item~\ref{claim:morestructure_item2} of Claim~\ref{claim:morestrucutre} we get $\{s,t\}\cap \{h_2,h_3\} = \emptyset$. Observe also that, by Itens~\ref{claim:morestructure_item3} and~\ref{claim:morestructure_item4}, if $s\in V(A_1)\cup \{h_1\}$, then $s$ must be equal to $h_1$ and, again by Item~\ref{claim:morestructure_item2}, we get that $A_1$ must be the multiedge $h_1h_2$. Because the analogous holds for $t$, we have the following possible cases:

	\begin{itemize}
		\item $s=h_1$ and $t=h_4$. By Lemma~\ref{lemma:simple_st} and Item~\ref{claim:morestructure_item3}, we have that $sh_2,h_2z_q,sz_0,z_0h_3,h_3t,z_qt$ all have multiplicity~1. Define: $a = \lambda(sz_0)$, $b = \lambda(z_0h_3)$, $c=\lambda(h_3t)$, $d=\lambda(sh_2)$, $e=\lambda(h_2z_q)$ and $f = \lambda(z_qt)$. 
		
		By Item~\ref{claim:morestructure_item1}, there exists a \tpath{s,t} using $z_0h_3$ that does not contain every vertex of degree~3. Observe that such path must be $P_1=(s,a,z_0,b,h_3,c,t)$. 
		Similarly, observe that the only \tpath{s,t} using $h_2z_q$ and not containing all vertices of degree~3 is $P_2 = (s,d,h_2,e,z_q,f,t)$. This is a contradiction as these paths are disjoint and $\{z_0,h_2\}$ is a \vcut{s,t} (i.e. $c(s,t) = 2$).
		
		\item $s=h_1$ and $t\in V(B_1)$. By Lemma~\ref{lem:edgecut}, we know that $V(B_1) = \{t\}$, by Lemma~\ref{lemma:simple_st}, that $sh_2$ and $h_2t$ have multiplicity~1, and by Item~\ref{claim:morestructure_item3} of Claim~\ref{claim:morestrucutre} that $z_qh_2$ has multiplicity~1. Therefore,  $c(s,t) = 2$ and $h_2$ is a common neighbor of $s$ and $t$ having degree 3 and such that every edge incident to $h_2$ has multiplicity~1. This contradicts Item~\ref{claim:morestructure_item4}.
	\end{itemize}

	\item $G$ is 2-crossed. We subdivide in a case analysis.
	
	\begin{itemize}
	    \item $c(s,t) = 3$: note that we get $d(s)\ge 3$ and $d(t)\geq 3$; hence $\{s,t\}\subseteq \{h_1,h_2,h_3,h_4\}$. Because these vertices are all symmetrical, suppose without loss of generality that $s = h_1$. 
	By Item~\ref{claim:morestructure_item2} and Lemma~\ref{lemma:simple_st}, we get $sh_2$ and $sh_4$ are multiedges of $G$ of multiplicity~1.  Additionally, since $st\notin E(G)$, we get that $t$ must be equal to $h_3$, which analogously as before gives us that $th_2$ and $th_4$  are multiedges of $G$ of multiplicity~1. 
	Now define the cut $S=\{h_2,z_0,h_4\}$, and observe that, for each $x\in S$, there is exactly only one $s,t$-path in $G - (S\setminus \{x\})$; denote such path by $P_{x}$. 
	Observe that the uniqueness of such paths imply that each of them must be a temporal path, as otherwise we could obtain a smaller cut. 
	We get a contradiction as $P_{z_0}$, $P_{h_2}$ and $P_{h_4}$ are 3 vertex disjoint \tpath{s,t}s.
	
	\item $c(s,t) =2$ and $\{s,t\} \subseteq \{h_1,h_2,h_3,h_4\}$. As in the previous case, we can suppose, without loss of generality, that $s=h_1$, which in turn gives us that $t=h_3$ and that $sh_2,sh_4,th_2,th_4$ are multiedges of $G$ of multiplicity~1. This and Item~\ref{claim:morestructure_item3} contradict Item~\ref{claim:morestructure_item4}. 
	
	\item $c(s,t) =2$ and $\{s,t\} \nsubseteq \{h_1,h_2,h_3,h_4\}$. First suppose that $s\in V(A_1)$. If $t \in B_{1}\cup B_2 \cup \{h_3,h_4\}$, by Lemmas~\ref{lemma:simple_st} and~\ref{lem:edgecut} and Item~\ref{claim:morestructure_item3} we get a situation that contradicts Item~\ref{claim:morestructure_item4} of Claim~\ref{claim:morestrucutre}. And because $st\notin E(G)$, we then get that $t\in V(A_2)$. 
	%
	Let $\alpha = \lambda(sh_1), \beta = \lambda(h_1z_0)$ and $\gamma = \lambda(f)$ where $f$ is an edge with endpoints $h_1 w$ for $w \in N(h_1)\setminus\{s,z_0\}$.
	First we prove that $\alpha \leq \gamma$. Suppose otherwise, then all \tpath{s,t}s using $f$ must also use $h_1 z_0$ and must start by using $s h_2$. If one of these paths uses $h_1 z_0$ before $f$, then we have $\beta \leq \gamma$. With the fact that $\alpha \leq \max \{\beta,\gamma\}$ (as otherwise no path uses $sh_1$), we get $\alpha \leq \gamma$, a contradiction. Therefore, every \tpath{s,t} that use $f$ also uses, in this order, the vertices $h_4,h_1,z_0$ and starts by $sh_2$. Let $P$ be such a path; we analyze the sequence of vertices of $P$ between the vertices $h_2$ and $h_4$. If $P$ goes to $h_3$ after visiting $h_2$, then as all $h_3,h_4$-paths that do not pass by $h_2$ must pass by $z_0$ or $t$, we know that $z_0$ is visited by $P$ before $h_4$, a contradiction. Therefore $P$ does not visit $h_3$ after $h_2$, which means that $P$ must start with $s,h_2,z_q$. 
	Now observe that, because $P$ must visit $h_4,h_1,z_0$ in this order, we can conclude that $P$ uses, in this order, $s,h_2,z_q,h_4,h_1,z_0,h_3,t$. 
	Therefore $P$ uses all vertices of degree~3, and since $P$ is a generic \tpath{s,t} using $f$, we know that this holds for all such paths, which contradicts Item~\ref{claim:morestructure_item1}. So, we have that $\alpha\le \gamma$, as we wanted to prove. Observe that this also gives us that there must exist a path passing by $f$ that visits $h_1$ before $h_4$, and since $f$ was an arbitrary edge with endpoints $h_1w$, we get that this holds for all such edges. 
	Note also that, since the graph is symmetric, this also holds for $h_2$, $h_3$ and related edges. The analogous arguments can be applied to the edges incident to $t$ to get that: $\lambda(th_3)\ge \lambda(e)$ for every $e$ with endpoints $h_3w_3$, where $w_3 \in N(h_3)\setminus\{t,z_0\}$; and $\lambda(th_4)\ge \lambda(e)$ for every $e$ with endpoints $h_4w_4$, where $w_4 \in N(h_4)\setminus\{t,z_0\}$.

	Now let $Q$ be any \tpath{h_1,h_4} and $Q'$ be any \tpath{h_2,h_3}. Since $\lambda(sh_1)\le \lambda(e)$ for every $e$ with endpoints $h_1w$, and $\lambda(th_4)\ge \lambda(e)$ for every $e$ with endpoints $h_4w_4$, we get that $P_1 = sQt$ is a \tpath{s,t}. Similarly we get that $P_2 = sQ't$ is a \tpath{s,t}, a contradiction as $P_1$ and $P_2$ are disjoint and $\{h_1,h_2\}$ is a \vcut{s,t} (i.e. $c(s,t) = 2$). 

	Finally, observe that the case where $s\in V(A_2)$ is analogous, as in this case we get $t\in V(A_1)$. Observe Figure~\ref{fig:ABtroca} to see that the cases $s\in V(B_1)$ or $s\in V(B_2)$ are also analogous. Additionally, if instead $t$ is within some of these subgraphs, then we can apply Lemma~\ref{lem:reverse}. 
	
	\begin{figure}[thb]
\begin{center}
\includegraphics[width=8cm]{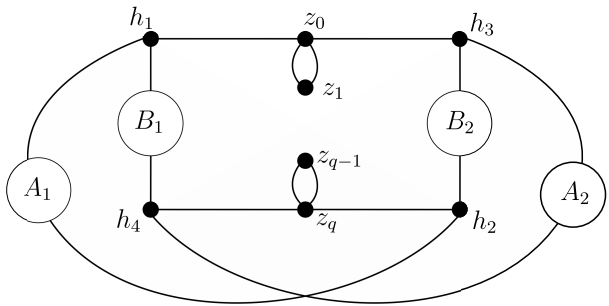}
\caption{A redrawing of $G$ that shows how the sets $A_1$ and $A_2$ are symmetric to $B_1$ and $B_2$.}
\label{fig:ABtroca}
\end{center}
\end{figure}

    \end{itemize}
    
\end{itemize}

\end{proof}



We are now ready to finish our proof.

\begin{proof}[Proof of sufficiency of Theorem~\ref{main}]
We want to prove that if $F\not\preceq_m G$ for every $F\in \{\mathcal{F}_1, \mathcal{F}_2, \mathcal{F}_3\}$, then $G$ is Mengerian.

	Suppose otherwise and let $(G,\lambda),s,t$ be a minimum counter-example. By Lemma~\ref{lemma:simple_st}, we know that $G$ is 2-connected; hence so is $U(G)$. 
	Observe that if $U(G)$ contains a subdivision of $\mathcal{F}_3$, then so also does $G$ and we are done. 
	Hence suppose that $U(G)$ has no ${\cal F}_3$ as topological minor. By Lemma~\ref{lem:Delta_geq_4} we get that $\Delta(G)\ge 4$, and by Lemma~\ref{lem:KKK.00} we get that there are $v,w\in V(G)$ and $d\ge 4$ for which $U(G)$ is $(v,w,d)$-decomposable.
	Observe that the components of $U(G)-\{v,w\}$ and of $G-\{v,w\}$ differ only on the multiplicity of the edges. 
	So let $G_{1},\cdots, G_p$ be the components of $G - \{v, w\}$, and for each $i \in [p]$, let $v_{i}$ (resp. $w_{i}$) be the vertex in $V(G_{i})\cap N(v)$ (resp. $V(G_{i})\cap N(w)$). Since $p\in \{d-1,d\}$ and $d\geq 4$, we have $p\geq 3$. 
	We now analyse the following cases.
	
	\begin{enumerate}
		
		\item\label{item:st=vw} $\{s,t\} = \{v,w\}$.  First suppose $s = v$ and $w = t$, and let $\mathcal{P}$ be a maximum set of vertex disjoint \tpath{s,t}s (i.e., $\lvert {\cal P}\rvert = p(s,t)$).  
		Since $s,t$ are non-adjacent, each path in ${\cal P}$ has at least one internal vertex. 
		So for each $P \in \mathcal{P}$, let $v_{i_P} = V(P)\cap \{v_{1},\cdots,v_{p}\}$ be the vertex in the neighborhood of $s$ used by $P$; also,  let $S = \{v_{i_P} ~\colon~ P \in \mathcal{P}\}$. Because the paths in ${\cal P}$ are internally disjoint, we know that $|S| = |{\cal P}|$. And since $\lvert {\cal P}\rvert = p(s,t) < c(s,t)$, we get that $S$ cannot be a \vcut{s,t}. Therefore, there exists a temporal path $P'$ in $(G-S,\lambda)$. Observe that $P'-\{s,t\}$ must be contained in $V(G_{j})$ with $v_j\notin S$, and that in this case $P'$ intersects the paths in $\mathcal{P}$ only in $\{s,t\}$. This contradicts the maximality of $\mathcal{P}$.

		For the case $s = t$ and $w=s$, we apply Lemma~\ref{lem:reverse} and the same argument as before to reach a contradiction.

		\item\label{item:stDistinctparts} There exist $i,j \in [p]$ with $i\neq j$ such that $s \in V(G_{i})$ and $t \in V(G_{j})$. Observe that both $\{vv_i,ww_i\}$ and $\{vv_j,ww_j\}$ satisfy Lemma~\ref{lem:edgecut}. Therefore, we have that $V(G_i) = \{s\}$ and $V(G_j) = \{t\}$. By Lemma~\ref{lemma:simple_st}, all multiedges incident to $s$ and $t$ have multiplicity~1. Hence, let $e_1$ be the edge with endpoints $sv$ and $e_2$ be the edge with endpoints $vt$. Because $w$ is not a \vcut{s,t}, there must be a \tpath{s,t} in $(G-w,\lambda)$. Additionally, because $G$ is $(v,w,d)$-decomposable, observe that such temporal path must be $(s,\lambda(e_1),v,\lambda(e_2),t)$. Applying a similar argument to $(G-v,\lambda)$, we get 2 disjoint \tpath{s,t}s in ${\cal G}$, a contradiction as $\{v,w\}$ is a \vcut{s,t}.

		\item $\{s,t\}\cap \{v,w\}$ is a unitary set. Suppose that $s = v$ and let $i\in [p]$ be such that $t\in V(G_i)$. In this case, again applying Lemma~\ref{lem:edgecut} to $\{sv_i,ww_i\}$, we get that $V(G_i) = \{t\}$, and hence $s$ and $t$ are adjacent, a contradiction. The other cases are clearly analogous.

		\item\label{item_stInGi}There exists $i\in [p]$ such that $\{s,t\}\subseteq V(G_i)$. This contradicts Lemma~\ref{lem:edgecut} applied to $\{vv_i,ww_i\}$. 
	\end{enumerate} 

\end{proof}

	\section{Mengerian Graphs Recognition}\label{sec:recog}
	
	In~\cite{GKMW.11}, they show that given two simple graphs $G$ and $H$, one can check if $H$ is a topological minor of $G$ in time $O(f(|V(H)|)|V(G)|^3)$. Thus, for a finite family of graphs ${\cal H} = \{H_{1},\ldots,H_{k}\}$, each $H_i$ of constant size, the problem of deciding whether $H\preceq G$ for some $H\in {\cal H}$ can be solved in polynomial time. If the same holds for m-topological minor, then Theorem~\ref{main} implies that we can recognize Mengerian graphs in polynomial time. We prove that this is indeed the case by giving an algorithm that makes use of the one presented in~\cite{GKMW.11}. Observe that we only need to recognize m-subdivisions of ${\cal F}_1$ and ${\cal F}_2$, since any m-subdivision of ${\cal F}_3$ is also a subdivision of ${\cal F}_3$ (i.e., we can recognize it by applying the algorithm in~\cite{GKMW.11}). The following lemma is the key for reaching polynomial time.
	
	\begin{lemma}\label{lem:chain}
		Let $G$ be a graph such that $F\preceq_m G$ for some $F\in \{{\cal F}_1,{\cal F}_2\}$ and $G$ has no $\mathcal{F}_3$ as topological minor. Then there is an m-subdivision $H\subseteq G$ of $F$ such that $d_{U(G)}(v)=2$ for every $v \in V(H)$ that is an internal vertex of the chain of $H$.
	\end{lemma}

\begin{proof}

We can suppose that $G$ is 2-connected as otherwise we can constrain ourselves to the 2-connected components of $G$. Let $F\preceq_m G$, where $F\in \{{\cal F}_1,{\cal F}_2\}$, and let $H\subseteq G$ be an m-subdivision of $F$ in $G$. Let $L = (z_{0},\ldots, z_{q})$ be the chain of $H$. We prove that if $d_{U(G)}(z_i)>2$ for some $i\in [q-1]$, then we can find an m-subdivision $H'\subseteq G$ of $F'\in \{{\cal F}_1,{\cal F}_2\}$ whose chain is smaller than the chain $L$. This means that if $H$ is an m-subdivision of a graph in $\{{\cal F}_1,{\cal F}_2\}$ with minimum chain size, then $d_{U(G)}(z_{i}) =2 $ for every $i \in [q-1]$ and the lemma follows. 

So let $z_i\in V(L)$ with $i\in [q-1]$ be such that $d_{U(G)}(z_i)> 2$, and let $v\in N(z_i)\setminus \{z_{i-1},z_{i+1}\}$. Because $U(G)$ is 2-connected as well, we get that $z_i$ cannot separate $v$ from $H$. Let $P$ be a path between $v$ and $H$ not passing through $z_i$, and let $u \in V(H)$ be the other endpoint of $P$. 
The proof consists of analyzing all possible choices of $F$, and where $u$ might lie in $H$.

\begin{itemize}
	\item $F = \mathcal{F}_{1}$. We use Proposition~\ref{prop:findF1} in $H$ to find subgraphs $C_1,C_2,J,L$ as indicated in Figure~\ref{fig:propf1}. 
	Let $w_{1} \in V(J)\cap V(C_{1})$ and $w_{2} \in V(J)\cap V(C_{2})$. 
	As the ordering of $L$ is arbitrary, we can suppose that $\{w_1,w_2\}\cap N(z_0) = \emptyset$. 
	
	First, consider that $u \in V(C_1)\setminus (N(z_0)\cup V(L))$, and observe Figure~\ref{fig:caso1-1} to follow the construction. 
	Let $Q$ be the $u,z_{i}$-path contained in $C_1$ not passing through $z_{i+1}$, and let $C'_1$ be the cycle formed by $P$ and $Q$. 
	Let $J'$ be either equal to $J$, in case $w_1\in V(C'_1)$, or be the $C'_1,C_2$-path defined by $J$ together with the $w_1,u$-path contained in $C_1 - V(L)$. 
	Finally, let $L'=(z_0,\ldots,z_i)$, and observe that $C'_1,C_2,J,L'$ satisfy the first two conditions in Proposition~\ref{prop:findF1}. The third condition also holds as $N_H(z_0)\cap V(C'_1) = N_H(z_0)\cap N(Q) = N_H(z_0)\cap N(C_1)$; thus we get a subdivision of ${\cal F}_1$ with a smaller chain, as desired.
	
	\begin{figure}
	\centering
	\includegraphics[width = 7cm]{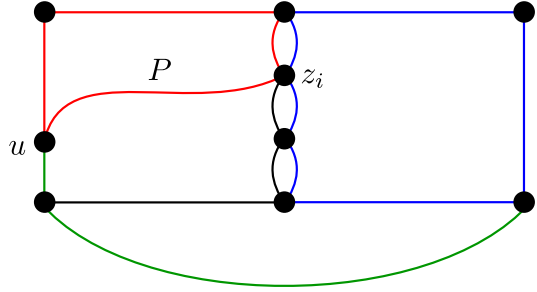}
	\caption{Case where $u\in V(C_1\setminus J)$. Cycle $C'_1$ is highlighted in red, cycle $C_2$ in blue, and path $J'$ in green.}\label{fig:caso1-1}
	\end{figure}
	
	Now suppose that $u\in (V(C_1)\cap N(z_0))\setminus V(L)$. Let $C'_1$ be the cycle defined by path $P$ plus the $u,z_i$-path contained in $C_1$ not passing through $z_0$, and let $L' = (z_i,\ldots,z_q)$. Again observe that  $C'_1,C_2,L',J$ satisfy the first two conditions of Proposition~\ref{prop:findF1}. For the third condition, observe that the neighbors of $z_i$ are contained in $V(P)\cup \{z_{i-1}\}$, which is disjoint from $\{w_1,w_2\}$.

	Consider now the case $u \in V(J)\setminus (V(C_1) \cup V(C_2))$. Observe Figure~\ref{fig:caso1-2} to follow the construction. Let $C'_1$ be the cycle formed by $P$ together with $uJw_1$ and the $w_1,z_i$-path contained in $C_1$ not passing by $z_q$. Also, let $L' = (z_0,\ldots, z_i)$, and $J'$ be equal to $uJw_2$. 
	One can again check that the conditions in Proposition~\ref{prop:findF1} hold.

		\begin{figure}[ht]
		\centering
		\includegraphics[width = 7cm]{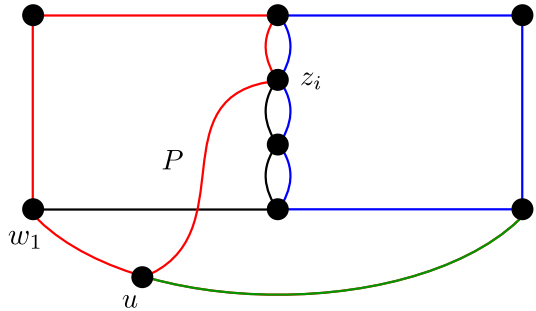}
		\caption{Case where $u\in V(J)\setminus (V(C_1) \cup V(C_2))$. Cycle $C'_1$ is highlighted in red, cycle $C_2$ in blue, and path $J'$ in green.}\label{fig:caso1-2}
	    \end{figure}
	
	Finally, suppose $u \in L$, i.e., $u=z_j$ for some $j\in\{0,\ldots,q\}$. Now we find a subdivision of $ \mathcal{F}_2$ with a smaller chain by using Proposition~\ref{prop:findF2}. 
	Observe Figure~\ref{fig:caso1-3} to follow the construction. Recall that $i>0$ and suppose first that $j>i$. Let $C'_{1}$ be the cycle formed by $P$ together with a $z_{i},z_{j}$-path contained in $L$. For each $\ell \in \{1,2\}$, let $J_{\ell}$ be the $w_{\ell},z_0$-path contained in $C_{\ell}$ that does not pass through $z_{q}$. Then define $C_{2}'$ as the cycle formed by the paths $J$, $J_{1}$ and $J_{2}$. 
	Finally, let $L'=(z_{0},\ldots,z_{i})$, and $J'$ be the $z_j,w_1$-path contained in $C_{1}$ that does not use vertex $z_{j-1}$. Recall that $P$ is a $u,v$-path, where $v\in N(z_i)\setminus \{z_{i-1},z_{i+1}\}$ to see that $C'_1$ is a cycle on at least~3 vertices. The other desired properties of $C'_2,L',J'$ can be seen to hold from the structure of $H$. Additionally, $L'$ is smaller than $L$, as desired.

    \begin{figure}[h]
	\centering
	\includegraphics[width = 7cm]{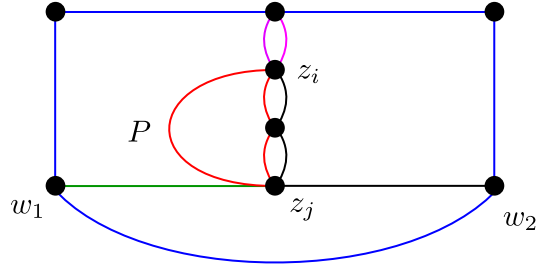}
	\caption{Case where $u = v_j\in V(L)$ for some $j\in \{i+1,\ldots,q\}$. Cycle $C'_1$ is highlighted in red, cycle $C'_2$ in blue, path $J'$ in green, and the new chain $L'$ in magenta.} \label{fig:caso1-3}
	\end{figure}
	
	For the case that $j<i$ we can apply a similar argument(see Figure~\ref{fig:caso1-3'}).
	
	\begin{figure}[h]
	\centering
	\includegraphics[width = 7cm]{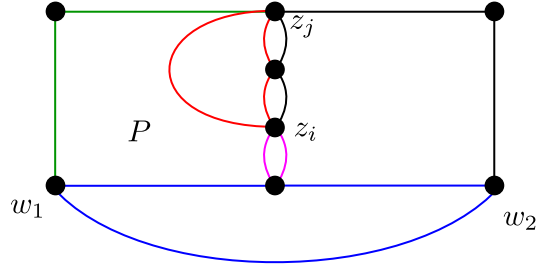}
	\caption{Case where $u = v_j\in V(L)$ for some $j\in \{0,\ldots,i-1\}$. Cycle $C'_1$ is highlighted in red, cycle $C'_2$ in blue, path $J'$ in green, and the new chain $L'$ in magenta.} \label{fig:caso1-3'}
	\end{figure}

	\item $F = \mathcal{F}_2$. Similarly as the last item, we start by using Proposition~\ref{prop:findF2}, to find $C_1,C_2$, linked by a chain $L$ and a path $J$, with $L$ and $J$ being disjoint (observe again Figure~\ref{fig:propf2} to see the related structures). Call $w_{1}$ the vertex in $V(J)\cap V(C_{1})$ and $w_{2}$ the vertex in $V(J)\cap V(C_{2})$. Additionally, let $V(C_{1})\cap V(L) = \{z_{0}\}$ and $V(C_{2})\cap V(L) = \{z_{q}\}$.

	Notice that the case where $u \in V(C_{1})\setminus V(L) $ is analogous to the case that $u \in V(C_{2}) \setminus V(L)$; so suppose $u\in V(C_1)\setminus V(L)$. 
	Let $C'_1$ be the cycle formed by $P$, together with the $u,z_0$-path contained in $C_1$ that also contains $w_1$, and a $z_0,z_i$-path contained in $L$. 
	Also, let $L' = (z_i,\ldots,z_q)$. 
	We have that $C'_1,C_2,L',J$ satisfy the properties in Proposition~\ref{prop:findF2} and $|V(L')|\leq |V(L)|$.
	
	Now suppose $u \in V(J)\setminus \{w_1, w_{2}\}$, and let $C'_1$ be the cycle formed by $P$, together with $w_1Ju$, any $w_1,z_0$-path contained in $C_1$, and a $z_0,z_i$-path contained in $L$. 
    Again, let $L' = (z_i,\ldots,z_q)$, and observe that $C_{1}',C_{2}, L', uJw_{2}$ satisfy Proposition~\ref{prop:findF2}.  
	
	Finally, suppose $u \in V(L)$, i.e., $u=z_j$ for some $j\in \{0,\ldots,q\}$. Observe that cycles $C_1$ and $C_2$ are symmetrical in the structure of $H$; hence, we can suppose, without loss of generality, that $j>i$. In this case, let $C_{2}'$ be the cycle formed by $P$ and a $z_i,z_j$-path contained in $L$; $J'$ be a $z_j,w_1$-path contained in $(L-\{z_0,\ldots,z_j\})\cup C_2\cup J$; and $L' = (z_0,\ldots,z_i)$. One can see that $C_1,C'_2,L',J'$ satisfy Proposition~\ref{prop:findF2}.  
	\end{itemize}
\end{proof}

\begin{proof}[Proof of Theorem~\ref{thm:recog}]
We need to prove that we can decide whether $G$ has one of the graphs in $\{\mathcal{F}_1, \mathcal{F}_2, \mathcal{F}_3\}$ as an m-topological minor in polynomial time where $n = \lvert V(G)\rvert$ and $m = \lvert E(G)\rvert$. 
First, we test if $G$ has $\mathcal{F}_3$ as topological minor, which can be done in time $O(n^3)$ according to the algorithm presented in~\cite{GKMW.11}. If the answer is positive, we are done as the concepts of m-topological minor and topological minor coincide because ${\cal F}_3$ is a simple graph. 

For a maximal chain $L$ such that the internal vertices of $L$ have degree~2 in $U(G)$, we identify all the vertices of $L$ into a single vertex, and test if the obtained graph, $G_L$, has a subdivision of $\mathcal{F}_3$. 
Observe that this operation applied to the chain of an m-subdivision of $\mathcal{F}_1$ or $\mathcal{F}_2$ would create a subdivision of $\mathcal{F}_3$. 
Therefore, if the answer is negative, then we can conclude that there is no m-subdivision of $\mathcal{F}_1$ or $\mathcal{F}_2$ whose chain is $L$. Otherwise, as $G$ itself has no subdivision of $\mathcal{F}_3$, while $G_L$ does, then we can apply Lemma~\ref{lem:helpcrossed} to either conclude that $G$ has $\mathcal{F}_1$ or $\mathcal{F}_2$ as topological m-minors, or that $G$ has a crossed structure. 
If the former occurs, we stop our algorithm. Otherwise, we apply the same procedure to another such chain. Finally, if the answer is negative for all such chains, then we can conclude from Lemma~\ref{lem:chain} that $G$ has no m-subdivisions of $\mathcal{F}_1$ or $\mathcal{F}_2$, and therefore must be a Mengerian graph. Since the total number of such chains is at most $m = \lvert E(G)\rvert$, we get that our algorithm runs in time $O(mn^3)$. 
\end{proof}
		
\bibliographystyle{plain}
\bibliography{references.bib}
	
\end{document}